\documentclass[11pt,oneside,final]{article}
\usepackage[a4paper, hmargin=2.5cm, vmargin=2.25cm]{geometry}
\usepackage[compact]{titlesec}
\usepackage{caption}
\usepackage{subcaption}
\usepackage{graphicx}
\usepackage{color, float}
\usepackage{listings}
\usepackage{setspace}
\usepackage{pdfpages}
\usepackage{natbib}
\bibpunct[, ]{(}{)}{,}{a}{}{,}

\usepackage{etoolbox}
\usepackage{rotating}
\usepackage{amsmath, amsthm, amsfonts, amssymb, amstext, mathabx}
\usepackage{appendix}
\usepackage{secdot}
\usepackage{accents}
\usepackage{makecell}
\usepackage[ruled,noend]{algorithm2e}
\usepackage{algorithmic}
\usepackage{array, cases}
\usepackage[colorlinks=true,citecolor=blue,linkcolor=blue]{hyperref}
\usepackage{multicol, multirow, epstopdf, booktabs, longtable}
\usepackage{bbm, enumitem, bm}
\usepackage{xr, xr-hyper}
\usepackage[capitalise,nameinlink]{cleveref}
\usepackage[acronym]{glossaries}

\makeatletter\def\thmheadbrackets#1#2#3{\thmname{#1}\thmnumber{
\@ifnotempty{#1}{ }\@upn{#2}}\thmnote{ {\the\thm@notefont[#3]}}}\makeatother
\newtheoremstyle{brakets}{}{}{\itshape}{}{\bfseries}{.}{ }{\thmheadbrackets{#1}{#2}{#3}}
\theoremstyle{brakets}

\newcommand\restr[2]{{\left.\kern-\nulldelimiterspace #1\right|_{#2}}}
\usepackage{url}
\makeatletter\g@addto@macro{\UrlBreaks}{\UrlOrds}\makeatother

\allowdisplaybreaks
\usepackage{lineno}
\usepackage{comment}
\usepackage{siunitx}
\title{Solving Integrated Periodic Railway Timetabling with Satisfiability Modulo Theories: A Scalable Approach to Routing and Vehicle Circulation}

\author{
    Florian Fuchs\textsuperscript{1},
    Bernardo Martin-Iradi\textsuperscript{1},
    Francesco Corman\textsuperscript{1,}\thanks{Corresponding author: francesco.corman@ivt.baug.ethz.ch}
}

\date{
    \textsuperscript{1}Institute for Transport Planning and Systems (IVT), ETH Zürich, 8092 Zürich, Switzerland
}

\begin{document}
\maketitle

\begin{spacing}{1.5}
\begin{abstract}
This paper introduces a novel approach for jointly solving the periodic \acrfull{ttp}, train routing, and \acrfull{vcp} through a unified optimization model. While these planning stages are traditionally addressed sequentially, their interdependencies often lead to suboptimal vehicle usage. We propose the \textit{VCR-PESP}, an integrated formulation that minimizes fleet size while ensuring feasible and infrastructure-compliant periodic timetables.

We present the first \acrfull{smt}-based method for the VCR-PESP to solve the resulting large-scale instances. Unlike the \acrfull{sat}, which requires time discretisation, \acrshort{smt} supports continuous time via difference constraints, eliminating the trade-off between temporal precision and encoding size. Our approach avoids rounding artifacts and scales effectively, outperforming both \acrshort{sat} and \acrfull{mip} models across non-trivial instances.

Using real-world data from the Swiss narrow-gauge operator \acrshort{rhb}, we conduct extensive experiments to assess the impact of time discretisation, vehicle circulation strategies, route flexibility, and planning integration. We show that discrete models inflate vehicle requirements and that fully integrated solutions substantially reduce fleet needs compared to sequential approaches. Our framework consistently delivers high-resolution solutions with tractable runtimes, even in large and complex networks.

By combining modeling accuracy with scalable solver technology, this work establishes \acrshort{smt} as a powerful tool for integrated railway planning. It demonstrates how relaxing discretisation and solving across planning layers enables more efficient and implementable timetables.
\end{abstract}


\end{spacing}
\clearpage
\begin{spacing}{1.5}
\section{Introduction}\label{sec:intro}

\begin{figure}[th!]
    \centering
    \includegraphics[width=0.6\textwidth]{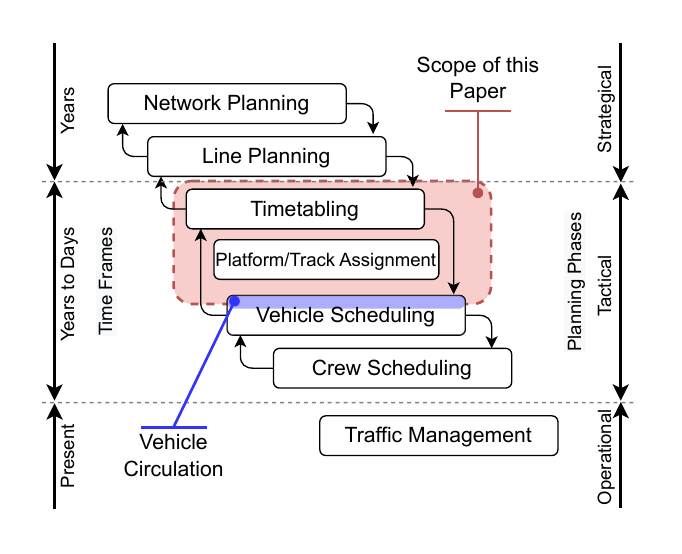}
    \caption{Illustration of railway planning stages, positioning vehicle circulation as an intermediate step between timetabling and vehicle scheduling.}
    \label{fig:planning_process}
\end{figure}

Efficient public transportation planning is crucial for urban mobility, and railway systems are one of its core modes of transport.
As depicted in \cref{fig:planning_process}, railway planning consists of strategic, tactical, and operational phases \citep{bussieck1997a, LUSBY20181}. At the strategic level, network and line planning define the overall service structure, determining routes, stops, and long-term infrastructure development. Tactical planning involves timetable design, vehicle assignment, and workforce allocation, whereas operational planning focuses on real-time traffic management to mitigate disruptions and maintain service reliability.
Vehicle circulation represents an intermediate planning stage between timetabling and detailed vehicle scheduling, where the focus is on determining the total number of vehicles needed rather than assigning specific rolling stock to individual trips.

Although planning these phases is traditionally handled sequentially, their interdependence suggests potential benefits from integrated optimization \citep{Schiewe2020IntegratedPlanning}. Joint optimization of timetabling and vehicle scheduling can reduce fleet requirements and improve resource utilization, highlighting the need for models that incorporate multiple planning aspects simultaneously.
However, solving these integrated problems using traditional \acrfull{mip} methods remains challenging due to computational complexity and achieving optimal solutions for large instances.

Recent studies \citep{Borndorfer2020AProblem,matos2020a,Gattermann2016IntegratingApproach} have shown that solution methods based on the \acrfull{sat} can outperform \acrshort{mip} solvers on planning problems such as the \acrfull{ttp}. Despite its advantages, \acrshort{sat} methods rely on the discretisation of time, which increases the size of the problem encoding as the time resolution becomes finer, posing scalability challenges.
\acrfull{smt} emerges as a promising alternative to address the scalability limitations. \acrshort{smt} extend \acrshort{sat} solvers by incorporating difference constraints, enabling the handling of continuous time without discretisation. This alternative allows for more scalable and precise solutions compared to traditional \acrshort{sat}-based methods.

This paper addresses all these limitations by, on one side, integrating timetabling, train routing, and vehicle circulation into a unified framework referred to as the Periodic Event Scheduling Problem with Vehicle Circulation and Routing (VCR-PESP). This breaks the sequential structure and focuses on minimizing the number of vehicles needed to operate a cyclic timetable. Conversely, we implement an \acrshort{smt}-based approach for solving the integrated problem at a large scale.

Key contributions of this research include:
\begin{itemize}
    \item \textbf{Integrated problem formulation:} We introduce the VCR-PESP, a unified model integrating periodic timetabling, infrastructure-aware routing, and vehicle circulation within a single optimization framework.
    
    \item \textbf{Scalable solution approach:} We develop a novel \acrshort{smt}-based model that supports continuous time and avoids discretisation. Compared to traditional \acrshort{sat} and \acrshort{mip} methods, it achieves improved scalability and runtime performance.
    
    \item \textbf{Real-world computational study:} We apply our method to real-world data from the Swiss narrow-gauge operator \acrshort{rhb}, demonstrating substantial fleet size reductions and efficient solver performance across different network and time resolutions.
    
    \item \textbf{Systematic evaluation of integration effects:} We quantify the impact of discretisation, routing flexibility, and vehicle sharing on vehicle requirements, illustrating the operational benefits of fully integrated planning approaches.
\end{itemize}

Our paper highlights the importance of integrating vehicle circulation into timetabling. It emphasizes the computational benefits of \acrshort{smt}, offering a framework for addressing the challenges of large-scale public transportation planning.

The paper is organised as follows: \Cref{s:litrev} provides a literature review covering timetabling, vehicle circulation, and train routing. \Cref{sec:methodology} details the methodology, focusing on integrating these components using \acrshort{smt}. \Cref{sec:results} presents and places the computational results in context. Finally, \Cref{sec:conclusion} summarises the findings and outlines directions for future research.

\section{Literature Review}\label{s:litrev}

This section reviews the stages of railway system planning, focusing on integrating periodic timetabling and vehicle circulation. We discuss recent advancements in optimization techniques and identify challenges that this paper seeks to address.

\subsection{Periodic Timetabling and Extensions, Including Routing}\label{subsec:lit_periodic_time_tabling}

As mentioned in Section \ref{sec:intro}, railway planning is divided into stages, with different stages handled sequentially due to their complex nature. In particular, the \acrshort{ttp} in the timetabling stage is a key planning problem, as it defines the arrival and departure times of all services and ensures that these do not conflict given the existing infrastructure and operational requirements.

A common approach for finding timetables for passenger trains in a network is to define a periodic timetable that repeats over a certain period. The \emph{Periodic Event Scheduling Problem} (PESP), introduced by \citet{Serafini1989AProblems}, is the foundation of periodic timetabling. It ensures that train schedules repeat in cycles while satisfying operational constraints. The \acrshort{ttp} framework has been widely applied in railway planning, balancing passenger demand, infrastructure capacity, and operational feasibility \citep{Peeters2003CyclicOptimization}. While PESP provides a structured optimization model, its reliance on modular constraints imposes limitations when incorporating additional degrees of freedom, such as train routing and infrastructure constraints.

Several extensions have sought to enhance PESP by improving timetable flexibility. Studies by \citet{Gattermann2016IntegratingApproach} and \citet{Robenek2016PassengerProblem} introduce passenger route adjustments, while \citet{Liebchen2004SymmetryTimetables} investigate timetable symmetry constraints. More recent work has focused on integrating track choice into periodic timetabling. \citet{Wust2019PeriodicInfrastructure} extend PESP to incorporate flexible train routing decisions, allowing trains to adapt to different infrastructure configurations. Similarly, \citet{MASING2023TimetablingTrack} explore routing adaptability in railway construction settings.

\citep{bortoletto_et_al2023-infra-aware} introduce the Infrastructure-aware PESP, where the model is formulated using explicit track-based constraints, and the Flexible Infrastructure Assignment PESP \citep{bortoletto_et_al2024-flex-infra}, generalizing infrastructure allocation across multiple configurations. 

Periodic timetabling can also be modeled using alternative models to the PESP.
Constraint-based formulations like the one presented in \citet{HEYDARStability2013171} incorporate multiple train types and explicitly minimize cycle time instead of relying on periodic constraints. \citet{robenek2016a} consider a model between cyclic and acyclic timetabling and display the interest of keeping regularity for a fixed set of lines instead of all lines to solve the problem on a network level. The model is solved using a simulated annealing method and verified on the Israeli railway network. \citet{MARTINIRADI2022511} use a time-space graph to formulate a macroscopic \acrshort{ttp} with train routing and solve it using a column-generation-based matheuristic.

\subsection{Integration of Vehicle Circulation into Timetabling}

Vehicle circulation is critical in railway operations \citep{borndrfer_et_al:OASIcs:2018:9721}, as it determines how rolling stock is allocated across scheduled services. Traditional approaches treat vehicle circulation as a separate problem, first establishing a timetable and then assigning vehicles, often resulting in inefficient fleet usage \citep{Goossens2006OnProblems}. More recent research has demonstrated that integrating vehicle circulation into \acrshort{ttp} can lead to reductions in fleet size while maintaining timetable feasibility \citep{lieshout2021VehiclePesp}. Optimization models can minimize turnaround times and dead-heading by considering rolling stock constraints at the timetabling stage, improving overall vehicle utilisation. 

Examples beyond railway include integrated models for bus operations that optimize fleet costs or passenger transfers, such as \citet{Ibarra_integrated}, \citet{FONSECA2018128}, and \citet{schmid2015integrated}. While these approaches explore valuable trade-offs between vehicle usage, robustness, and transfers, they typically assume fixed routing or ignore infrastructure constraints, making them less applicable to railway models.

\acrshort{mip} formulations are widely used for periodic timetabling and vehicle circulation problems due to their flexibility in constraint modeling \citep{Goerigk2017AnProblem, Herrigel2018PeriodicModel}. However, due to their inherent complexity, \acrshort{mip} models become computationally intractable for large-scale instances.

Even without additional constraints such as routing or rolling stock optimization, PESP is known to be NP-hard \citep{Peeters2003CyclicOptimization}. This computational intractability has motivated the development of new solution techniques, relying on \acrfull{sat}.

\subsection{\acrshort{sat}-based Methods in Transport Planning}

Recent studies have shown that solving the \acrshort{ttp} without \acrshort{vcp} using \acrfull{sat}-based solvers outperforms commercial ones for \acrshort{mip}s, especially when the objective is to find feasible solutions \citep{grossmann2012a, Gromann2016SatisfiabilityProblems, Kummling2015AComputation, FuchsTrivellaCorman2022}. \citet{Borndorfer2020AProblem} propose a concurrent solver for the \acrshort{ttp} that integrates \acrshort{sat}, \acrshort{mip}, and domain-specific heuristics, demonstrating the benefit of combining logical and mathematical programming approaches for solving PESP instances.

\acrshort{sat}-based methods encode scheduling constraints as logical formulas and leverage efficient SAT solvers to determine feasible solutions. However, \acrshort{sat} formulations rely on time discretisation, which increases encoding size and computational effort as resolution granularity improves \citep{Gromann2011PolynomialSATb}. While this can be mitigated through coarse time steps, such approximations may lead to suboptimal or infeasible solutions in high-resolution timetables.

Previous extensions of SAT include the work by \citet{Gattermann2016IntegratingApproach}, who augment SAT-based timetabling by including passenger routing preferences in the form of soft constraints, resulting in a MaxSAT formulation. \citet{matos2020a} further enhance this approach by introducing reinforcement learning to guide the search in MaxSAT, improving performance on public benchmark instances. While these methods demonstrate the adaptability of logic-based frameworks to include user-centric objectives, they do not yet support integration with routing and vehicle constraints, which introduces additional modeling and computational challenges.

These developments illustrate a shift from classical SAT solving toward hybrid and learning-based strategies, highlighting the potential of logic-oriented approaches to model passenger-centric objectives. However, such approaches have not yet been extended to integrate routing or vehicle constraints, which introduces additional structural complexity.

\subsection{Research Gap}

While prior research has addressed periodic timetabling extensions, routing flexibility, and vehicle circulation separately, their joint optimization remains an open challenge. Most existing models treat routing independently of rolling stock constraints or optimize vehicle circulation, assuming fixed train paths \citep{Caimi2017ModelsPractice}.

Moreover, although \acrshort{mip} and \acrshort{sat}-based methods have successfully solved isolated components of the problem, they struggle with scalability when considering large-scale instances with high-resolution schedules. In particular, \acrshort{sat}-based models require time discretisation, which increases encoding size and computational complexity as resolution improves \citep{Gromann2011PolynomialSATb}. While coarse time steps can reduce problem size, they may lead to suboptimal or infeasible solutions.

 This paper formulates an integrated optimization model that addresses these challenges. The model jointly considers \acrshort{ttp}, routing flexibility, and vehicle circulation at a mesoscopic resolution—offering more detail than macroscopic models while remaining computationally tractable. We leverage Satisfiability Modulo Theories (\acrshort{smt}) as a scalable solution framework. \acrshort{smt} extends \acrshort{sat} by incorporating difference-logic constraints, allowing continuous-time formulations without discretisation \citep{armando2004sat, LEUTWILER2022525}. This extension enables efficient solving of integrated problems while maintaining high temporal precision.

\section{Methodology}\label{sec:methodology}

This section presents the methodology for modeling and optimizing train routes, event timings, train sequencing, and vehicle transitions within a periodic timetable. We introduce the \acrfull{ean}, the foundational structure for modeling events, activities, and interdependencies. The \acrshort{ean} provides a structured representation of the timetabling problem, linking routing and scheduling decisions.

Building on the \acrshort{ean}, we describe three distinct approaches to solving the problem:
a MIP formulation, an SMT encoding, and a SAT-based version. Each approach builds on the \acrshort{ean} with selectable activities, enabling us to solve and compare the three.

\subsection{Introducing the \acrfull{ean}}\label{ss:ean}

The \acrshort{ean} is a standard data structure to model periodic timetabling problems. It consists of a set of nodes $\mathcal{E}$ representing \emph{events} (e.g., train arrivals or departures) and a set of edges $\mathcal{A}$ representing \emph{activities} that capture constraints between these events \citep{Liebchen2004TheBeyond}. While activities are traditionally binding and consistently enforced, enabling routing flexibility or optional passenger and vehicle transfers requires defining some activities as \emph{selectable}—their constraints apply only if specific conditions are met. This concept supports modeling infrastructure-dependent routing \citep{FuchsTrivellaCorman2022} and conditional transfers \citep{Kroon2014FlexibleTimetabling}, and is essential for integrated optimization across planning stages.

\begin{figure}[H]%
\begin{center}
\includegraphics[width=0.6\textwidth]{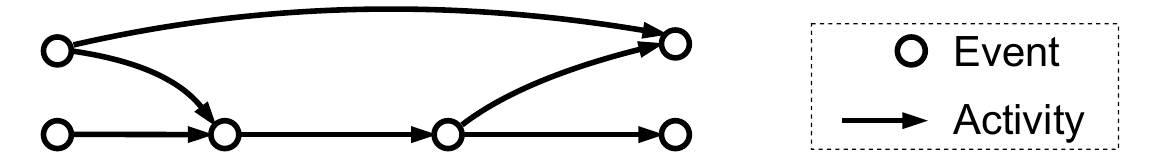}
\caption{Example \acrshort{ean} with six events and six activities.}
\label{fig:ean_example}
\end{center}
\end{figure}

The \acrshort{ean} aims to assign a timestamp $t_e \in [0, T)$ to each event $e \in \mathcal{E}$, where $T$ is the period of the timetable. For each activity $a \in \mathcal{A}$, there is a duration $\delta_a$ bounded by a lower limit $\delta_a^{\text{min}}$ and an upper limit $\delta_a^{\text{max}}$, such that $\delta_a^{\text{min}} \leq \delta_a \leq \delta_a^{\text{max}}$. The duration $\delta_a$ is derived from the scheduled times of the two connected events. If $i$ and $j$ denote the origin and destination events of activity $a=(i,j)$, then the duration is calculated as:
\begin{align}
\delta_a = t_{j} - t_{i} + k_a \cdot T, \quad \forall a=(i,j) \in \mathcal{A}
\end{align}
where $k_a \in \mathbb{Z}$ adjusts for cases where $t_{i} > t_{j}$ due to train schedules being longer than the cycle period.

The \acrshort{ean} includes three types of events: \emph{arrival}, where trains arrive at a stopping stop, \emph{departure}, where trains depart from a stopping stop, and \emph{passing}, where a train traverses a given location. It also supports six distinct activity types:
\begin{enumerate}
    \item \textit{Trip activities} for train movements between stations,
    \item \textit{Dwell activities} for stopping at stations,
    \item \textit{Headway activities} to ensure safe and feasible separation of trains using shared infrastructure,
    \item \textit{Regularity activities} to enforce even spacing for high-frequency services,
    \item \textit{Commercial activities} to meet service requirements such as maximum time between commercial stops,
    \item \textit{Vehicle-transfer activities} to model vehicle handovers between terminating and originating trains at termini \citep{Kroon2014FlexibleTimetabling, lieshout2021VehiclePesp}.
\end{enumerate}
These activities are explained in detail in the remainder of this section.

\subsection{Including Train Routing}\label{sec:tfn}

\begin{figure}[H]%
\begin{center}
\includegraphics[width=0.5\textwidth]{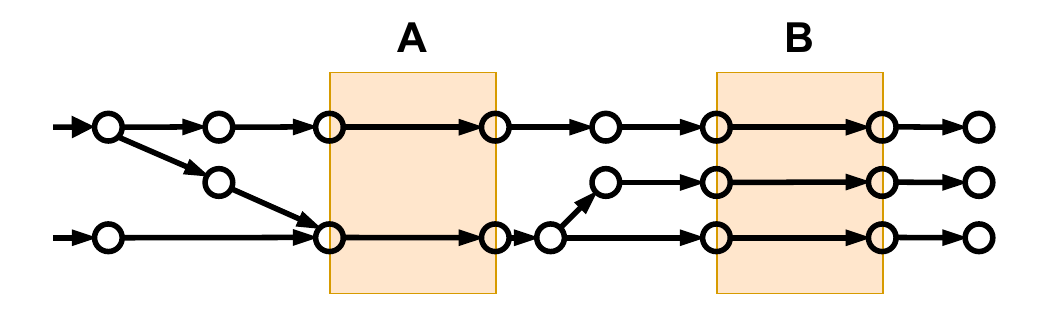}
\caption{Commercial and itinerary constraints for an example train visiting two stops, A and B, modeled as activities in the \acrshort{ean}.}
\label{fig:cc_ean_example}
\end{center}
\end{figure}

To optimize train routes, we extend the \acrshort{ean} to include multiple routing options per train service, following the approach of \citet{FuchsTrivellaCorman2022}. Routes consist of \textit{dwell} and \textit{trip} activities, which are modeled within the \acrshort{ean}. An example of itinerary activities for a single train is depicted in \cref{fig:cc_ean_example}.

\begin{figure}[H]%
\begin{center}
\includegraphics[width=0.5\textwidth, trim=0 0 0 0, clip]{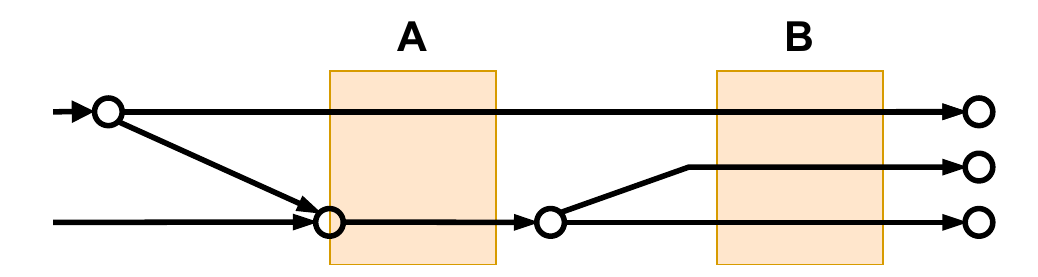}
\caption{A  \acrshort{tfn} example for a train service with five possible paths.}
\label{fig:TFN}
\end{center}
\end{figure}

To construct this extended \acrshort{ean}, we begin with a line plan specifying the stations each train serves. Each train's infrastructure and routing options are derived from this station sequence, allowing the \acrshort{ean} to represent the train's feasible movements through the network. The resulting network includes one weakly connected component for each train, which consists of only \textit{dwell} and \textit{trip} activities. These are organized into a \acrfull{tfn}, a graph with node set $\mathcal{V}$ and arc set $\mathcal{W}$ representing itinerary options. An example \acrshort{tfn} is shown in \cref{fig:TFN}.

The sets $\mathcal{V}^{\text{Source}}$ and $\mathcal{V}^{\text{Sink}}$ denote artificial source and sink nodes for the set of train services to schedule, and are a subset of $\mathcal{V}$.
The \acrshort{tfn} (as well as the \acrshort{ean}) is a directed and acyclic graph, meaning that a path from source to sink will define a valid route for a train service.
To compute the chosen route, we define binary variables
$x_v \in \{0,1\}$ to indicate whether a node $v \in \mathcal{V}$ is visited, and $x_w \in \{0,1\}$ to indicate whether an arc $w \in \mathcal{W}$ is selected. 
As the \acrshort{tfn} is a reduced version of the \acrshort{ean}, each event $e \in \mathcal{E}$ and activity $a \in \mathcal{A}$ can be mapped to the related nodes and links in the \acrshort{tfn}. Once the train routes are selected, the \acrshort{ean} is completed by incorporating all remaining activities, such as \textit{headway} and \textit{vehicle-transfer} activities. This completion ensures that the final model respects operational and commercial constraints while preventing conflicts.

\subsection{Adding Vehicle Circulation}\label{ss:vehicle_circulation}

Optimized vehicle circulation ensures the efficient use of rolling stock. To model vehicle movements effectively, we extend the \acrshort{ean} to include vehicle circulation links at termini. These links represent the transfer of vehicles between train services at their origin and destination stations. By modeling these links as activities within the \acrshort{ean}, we can seamlessly apply the same framework for routing trains to vehicle circulation.

\begin{figure}[H]
    \centering
    \begin{subfigure}[t]{0.4\textwidth}
        \centering
        \includegraphics[width=\textwidth,trim={0 0 20cm 0},clip]{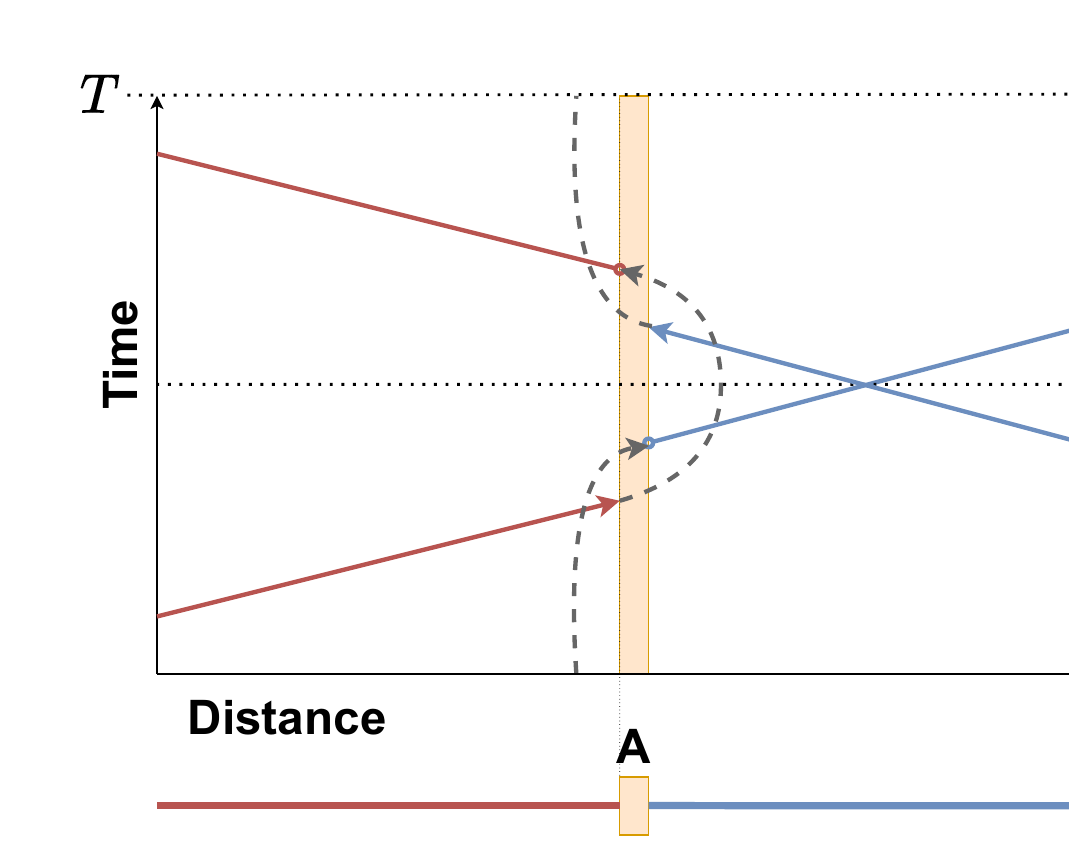}
        \caption{Vehicle circulation restricted to individual lines.}
        \label{fig:vehicle_per_line}
    \end{subfigure}
    \hspace{0.1\textwidth} 
    \begin{subfigure}[t]{0.4\textwidth}
        \centering
         \includegraphics[width=\textwidth,trim={20cm 0 0 0},clip]{figures/VehicleTransfers.drawio.pdf}
        \caption{Vehicles shared across lines.}
        \label{fig:vehicle_shared}
    \end{subfigure}
    \caption{Illustration of vehicle circulation strategies using two lines (i.e., red and blue): per-line circulation (a) and shared vehicle circulation (b).}
    \label{fig:vehicle_circulation_modes}
\end{figure}

\Cref{fig:vehicle_circulation_modes} illustrates the two primary vehicle circulation strategies. In the \textit{per-line} strategy (\cref{fig:vehicle_per_line}), vehicles are restricted to operating within the same line, effectively isolating vehicle pools for each train line. In the \textit{shared circulation} strategy (\cref{fig:vehicle_shared}), vehicles can transfer freely between lines, enabling a more flexible allocation of rolling stock. The choice of strategy may have significant implications for vehicle requirements and operational flexibility, as suggested by the example in \cref{fig:vehicle_circulation_modes} and demonstrated in our experimental results.

For each terminating train service, the model considers all possible links that transfer vehicles to originating train services. These links are represented as binary variables \( x_u \in \{0, 1\} \), where \( u \in \mathcal{U} \) denotes a vehicle circulation link. 
To ensure a feasible matching between train arrivals and departures, we impose the following constraints:
\begin{align}
    \sum_{u \in \mathcal{U}^{\text{Out}}_e} x_u &= 1, \quad \forall e \in \mathcal{E}^{\text{Start}}, \label{eq:vehicle_departure_circulation_ex} \\
    \sum_{u \in \mathcal{U}^{\text{In}}_e} x_u &= 1, \quad \forall e \in \mathcal{E}^{\text{End}}, \label{eq:vehicle_arrival_circulation_ex}
\end{align}
where the events of all originating and terminating train services are denoted by the sets $\mathcal{E}^{\text{Start}}$ and $\mathcal{E}^{\text{End}}$, respectively, \( \mathcal{U}^{\text{Out}}_e \) is the set of outgoing circulation links associated with event \( e \), and \( \mathcal{U}^{\text{In}}_e \) is the set of incoming circulation links. These constraints ensure that each arrival event at a terminus matches precisely one departure event, forming a valid vehicle transfer. Due to the cyclic characteristic of the timetable, we do not need to account for the start and end of vehicle operations.

\subsection{Defining Selectable Activities}\label{ss:selectable}

We follow the concept of \emph{selectable} and \emph{non-selectable} activities introduced by \citet{FuchsTrivellaCorman2022}, and we extend this approach by treating all activities $a \in \mathcal{A}$ as \emph{selectable}. A \emph{selectable} activity is one whose time constraints $\delta_a \in [\delta_a^{\text{min}}, \delta_a^{\text{max}}]$ are enforced only if associated routing or circulation decisions activate the activity. Otherwise, the activity remains inactive and its duration unconstrained~(i.e., $t_a \in \mathbb{R}$). As mentioned in Section \ref{sec:tfn}, each activity $a \in \mathcal{A}$ has associated elements in the \acrshort{tfn}, which group together a set of itinerary nodes $v \in \mathcal{V}_a \subseteq \mathcal{V}$, links $w \in \mathcal{W}_a \subseteq \mathcal{W}$, or vehicle transfers $u \in \mathcal{U}_a \subseteq \mathcal{U}$, and the activation of the activity depends on the utilization of such elements.

To formalize this relationship, for each activity $a \in \mathcal{A}$, we define its relevant components $\Omega_a$ as the union of the sets $\mathcal{V}_a$, $\mathcal{W}_a$, and $\mathcal{U}_a$:
\begin{align}
    \Omega_a = \mathcal{V}_a \cup \mathcal{W}_a \cup \mathcal{U}_a. \label{eq:union_of_activations}
\end{align}
The sets $\mathcal{V}_a$, $\mathcal{W}_a$, and $\mathcal{U}_a$ are activity-specific. For example:
\begin{itemize}
    \item A \textit{trip} activity depends solely on the corresponding infrastructure link $w \in \mathcal{W}$, with $|\mathcal{U}_a| = 0$, $|\mathcal{V}_a| = 0$, and $|\mathcal{W}_a| = 1$.
    \item A \textit{headway} activity depends on two infrastructure nodes $v \in \mathcal{V}$ traversed by the trains, with $|\mathcal{U}_a| = 0$, $|\mathcal{V}_a| = 2$, and $|\mathcal{W}_a| = 0$.
    \item A \textit{vehicle-transfer} activity depends on a single vehicle transfer, with $|\mathcal{U}_a| = 1$, $|\mathcal{V}_a| = 0$, and $|\mathcal{W}_a| = 0$.
\end{itemize}

The duration $\delta_a$ of an activity $a \in \mathcal{A}$, is by default constrained by the time bounds $\delta_a^{\text{min}}$ and $\delta_a^{\text{max}}$ if all associated components $\omega \in \Omega_a$ are active. If at least one component is not selected ($x_\omega = 0$), these bounds are relaxed using large constants $B_a^{\text{min}}$ and $B_a^{\text{max}}$:
\begin{align}
    \delta_a &\geq \delta_a^{\text{min}} - (|\Omega_a| - \sum_{\omega \in \Omega_a} x_\omega) \cdot B_a^{\text{min}}, \quad \forall a \in \mathcal{A}, \label{eq:selectable_activities_lb} \\
    \delta_a &\leq \delta_a^{\text{max}} + (|\Omega_a| - \sum_{\omega \in \Omega_a} x_\omega) \cdot B_a^{\text{max}}, \quad \forall a \in \mathcal{A}. \label{eq:selectable_activities_ub}
\end{align}
These equations ensure that the bounds on $\delta_a$ are only enforced when the activity is active, while inactive activities are effectively unconstrained.

This unified treatment of selectable activities eliminates the need to explicitly distinguish between \emph{selectable} and \emph{non-selectable} activities, as all activities are inherently treated as selectable. For example, independent of routing decisions, \textit{commercial} activities are modeled with $|\Omega_a| = 0$. Consequently, their time bounds must always be respected.
For notation simplicity, we unify the definition of our selection variable $x_\omega$ for each element $\omega \in \Omega_a$ and activity $a \in \mathcal{A}$.

By linking the activation of activities to the \acrshort{tfn}, we effectively reduce redundancy in the model. For instance, headway activities are automatically deactivated if only one train uses the relevant infrastructure, and redundant vehicle-transfer options are excluded based on the chosen routes. This approach prevents over-constraining the problem, ensuring feasibility while maintaining flexibility for optimization.

\subsection{\acrshort{mip} Formulation}\label{ss:mip_form}

As a final step before formulating the \acrshort{mip}, we define the objective function, which minimizes the number of vehicles required to operate a periodic timetable. In the VCR-PESP, each train service must be part of a feasible vehicle sequence that loops back periodically, forming a closed path along a cycle in the \acrshort{ean}. These vehicle cycles consist of alternating commercial activities (i.e., scheduled train trips) and vehicle-transfer activities (i.e., transitions at terminus stations).

Let \( \mathcal{S} \subset \mathcal{A} \) be the set of commercial activities and \( \mathcal{U} \subset \mathcal{A} \) the set of selected vehicle-transfer activities. When a valid timetable is constructed, each vehicle must follow a path through a sequence of activities in \( \mathcal{S} \cup \mathcal{U} \) that forms a cycle of total duration equal to an integer multiple of the period \( T \). This property reflects the cyclic nature of the timetable, where vehicles repeat the same circulation pattern every period.

\begin{figure}[H]%
\begin{center}
\includegraphics[width=0.69\textwidth]{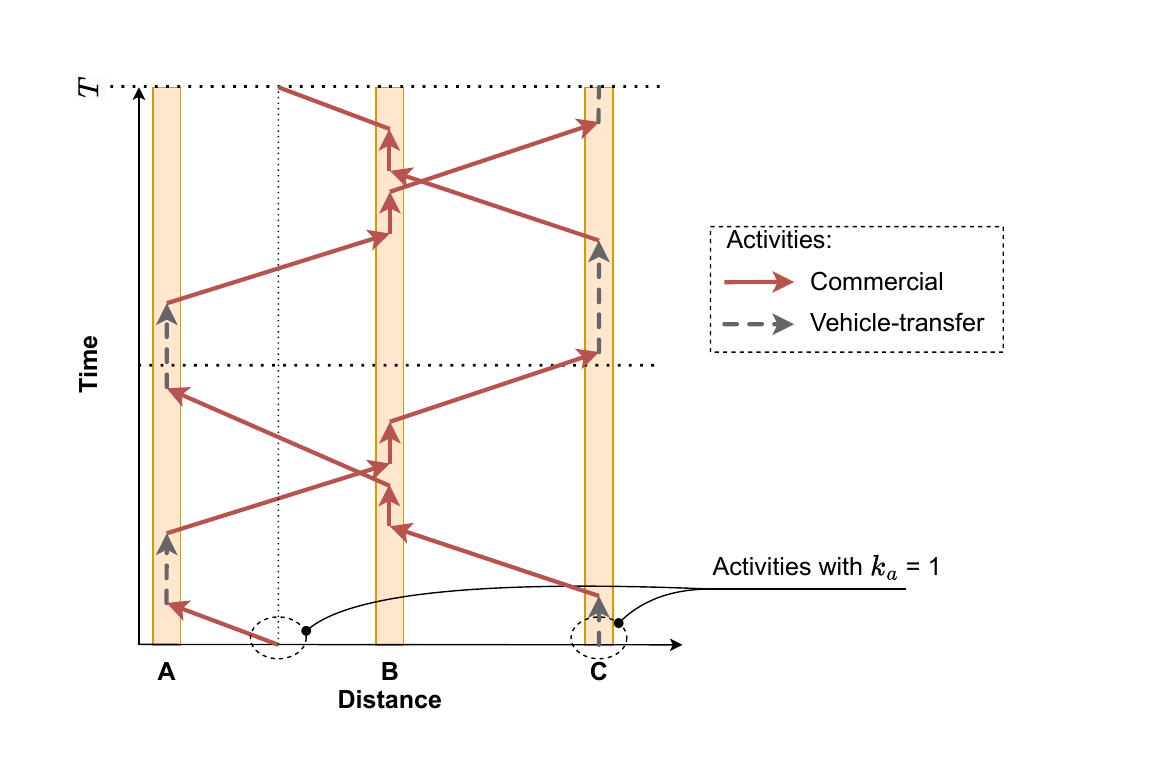}
\caption{Example illustration of two cycles in which vehicles circulate. Each cycle has a total length of \( T \). Since there are two such cycles, the total vehicle requirement is 2.}
\label{fig:vehicle_cycle}
\end{center}
\end{figure}

An illustrative example is shown in \cref{fig:vehicle_cycle}, highlighting two such vehicle cycles. Each cycle consists of alternating activities with duration \( \delta_a = t_j - t_i + k_a T \). Because each event appears exactly once as a predecessor and once as a successor in the cycle, all \( t_j - t_i \) terms cancel when summing over the cycle. The total duration of the cycle reduces to \( \sum_{a \in \text{cycle}} k_a \cdot T \), and the number of vehicles needed to operate the cycle is thus equal to \( \sum_{a \in \text{cycle}} k_a \).

\noindent
Consequently, the total number of vehicles required across all cycles in the solution is given by:
\[
\sum_{a \in \mathcal{S} \cup \mathcal{U}} k_a
\]
This sum reflects the objective function of the VCR-PESP model and is minimized to obtain the most efficient vehicle circulation.

\begin{table}[H]
\centering
\caption{Consolidated Variable and Set Definitions for VCR-PESP}
\label{tab:consolidated_variables_and_sets}
{\footnotesize \resizebox{\textwidth}{!}{\begin{tabular}{lll}
\hline
Symbol & Type & Definition \\ \hline
$\mathcal{V}$ & Set & Set of itinerary nodes in the \acrshort{tfn}. \\
$\mathcal{V}^{\text{Source}}$, $\mathcal{V}^{\text{Sink}}$ & Set & Set of artificial source and sink nodes for train services in the \acrshort{tfn}. \\

$\mathcal{W}$ & Set & Set of itinerary links in the \acrshort{tfn} \\
$\mathcal{U}$ & Set & Set of vehicle circulation links. \\
$\mathcal{U}^{\text{In}}_e, \mathcal{U}^{\text{Out}}_e,$ & Set & Set of outgoing (resp. incoming) circulation links associated with event \( e \in \mathcal{E}^{\text{End}} \) (resp. \( e \in \mathcal{E}^{\text{Start}} \)). \\
$\mathcal{S}$ & Set & Set of commercial activities of train services. \\
$\mathcal{E}$ & Set & Set of events (i.e., arrival, departure, passing) in the \acrshort{ean}. \\
$\mathcal{E}^{\text{Start}}$, $\mathcal{E}^{\text{End}}$ & Set & Set of originating and terminating train service events. \\
$\mathcal{A}$ & Set & Set of activities (e.g., trip, dwell, headway, regularity, commercial, vehicle-transfer) in the \acrshort{ean}. \\
$\mathcal{V}_a$, $\mathcal{W}_a, \mathcal{U}_a$ & Set & Sets of itinerary nodes and links, and vehicle circulation links associated with activity $a \in \mathcal{A}$. \\
$\Omega_a$ & Set & Union of sets $\mathcal{V}_a, \mathcal{W}_a, \mathcal{U}_a$ for activity $a \in \mathcal{A}$. \\ 
$\alpha^+(v), \alpha^-(v)$ & Set & Set of outgoing (resp. incoming) nodes from (resp. to) $v \in \mathcal{V}$ directly connected by a link in $\mathcal{W}$. \\ \hline
$T$ & Parameter & Period duration. \\
$\delta_{a}^{\text{min}}$ & Parameter & Lower bound for duration of activity $a \in \mathcal{A}$. \\
$\delta_{a}^{\text{max}}$ & Parameter & Upper bound for duration of activity $a \in \mathcal{A}$. \\
$B_{a}^{\min}$, $B_{a}^{\max}$ & Parameter & Bounds for relaxed duration of non-selected activity $a \in \mathcal{A}$. \\ \hline
$x_\omega$ & Variable & Binary variable for selection of element $\omega \in \Omega_a$ of activity $a \in \mathcal{A}$. \\
$t_e$ & Variable & Timestamp of event $e \in \mathcal{E}$ \\
$\delta_a$ & Variable& Duration of activity $a \in \mathcal{A}$ \\
$k_a$ & Variable & Period adjustment for activity $a \in \mathcal{A}$ \\
\hline
    \end{tabular}}}
\end{table}
We present in \eqref{pesp:formulation} the formulation of the VCR-PESP. The notation is summarized in \cref{tab:consolidated_variables_and_sets}

\begin{subequations}
\label{pesp:formulation}
\begin{align}
\min & \sum_{a \in \mathcal{S} \cup \mathcal{U}} k_a \label{eq:pesp_objective}\\ 
\text{subject to:\,\,}  &  \delta_a \geq \delta_{a}^{\min} - (|\mathcal{V}_a| - \sum_{v \in \mathcal{V}_a} x_v + |\mathcal{W}_a| - \sum_{w \in \mathcal{W}_a} x_w ) \cdot B_{a}^{\min}, & \forall{a} \in \mathcal{A}, \label{eq:pesp_lb_activation} \\  
&   \delta_a \leq \delta_{a}^{\max} +  (|\mathcal{V}_a| - \sum_{v \in \mathcal{V}_a} x_v + |\mathcal{W}_a| - \sum_{w \in \mathcal{W}_a} x_w ) \cdot B_{a}^{\max}  & \forall{a} \in \mathcal{A}, \label{eq:pesp_ub_activation} \\
& \delta_a = (t_{j} - t_{i}) + k_{a} \cdot T,  & \forall a =(i,j) \in \mathcal{A},\label{eq:pesp_activity_span}\\
&  \sum_{w \in \alpha^+(v)} x_w = \sum_{w \in \alpha^-(v)} x_w = x_v, & \forall v \in \mathcal{V},\label{eq:flow_conservation} \\
& \sum_{w \in \alpha^+(v)} x_w = 1,&  \forall v \in \mathcal{V}^{\text{Source}} \label{eqs:demand_one_train} \\
& \sum_{u \in \mathcal{U}^{\text{Out}}_e} x_{u} = 1, & \forall e \in \mathcal{E}^{\text{Start}}, \label{eq:vehicle_departure_circulation} \\
& \sum_{u \in \mathcal{U}^{\text{In}}_e} x_{u} = 1, & \forall e \in \mathcal{E}^{\text{End}}, \label{eq:vehicle_arrival_circulation} \\
\text{variables:\,\,} &  t_e \in [0, T), & \forall e \in \mathcal{E}, \label{eq:event_time} \\
&  \delta_a \geq 0, & \forall a \in \mathcal{A} , \label{eq:variables}\\
&  x_{\omega} \in \{0,1\}, & \forall \omega \in \Omega_a, \forall a \in \mathcal{A}, \label{eq:vehicle_circulation_variables} \\
&  k_a \in \mathbb{Z}, & \forall a \in \mathcal{A} \label{eq:variables3}
\end{align}
\end{subequations}

The objective function in \eqref{eq:pesp_objective} minimizes the vehicle count, as suggested by \citet{lieshout2021VehiclePesp}. Constraints \eqref{eq:pesp_lb_activation} and \eqref{eq:pesp_ub_activation} ensure both the correct activation of elements and the duration window of activities. The PESP constraint, given in \eqref{eq:pesp_activity_span}, establishes a direct relationship between the duration of an activity and the timing of its associated events. 
Constraint \eqref{eq:flow_conservation} guarantees flow conservation in the \acrshort{tfn}, and together with Constraint \eqref{eqs:demand_one_train}, it enforces that exactly one route is chosen for each train service, thus requiring a unique path for each train within the timetable. Vehicle circulation is managed through constraints \eqref{eq:vehicle_departure_circulation} and \eqref{eq:vehicle_arrival_circulation}. These ensure that for each event that marks the start or end of a train's journey, precisely one vehicle circulation activity is chosen, as proposed by \citet{lieshout2021VehiclePesp}. The nature of the decision variables used in the model is defined in \eqref{eq:event_time}, \eqref{eq:variables}, \eqref{eq:vehicle_circulation_variables}, and \eqref{eq:variables3}. 

\subsection{Transformation to \acrshort{smt}}

The VCR-PESP can quickly become intractable for state-of-the-art commercial solvers. To solve \eqref{pesp:formulation}, we present an alternative approach using \acrfull{smt}, which integrates \acrfull{sat} with difference constraints. This approach allows us to encode the model formulation~\eqref{pesp:formulation} into a combination of Boolean formulas and arithmetic constraints. Specifically, the Boolean formula is expressed in its \acrfull{cnf} and consists of Boolean variables \( q \in \mathcal{Q} \), which form clauses composed of literals. Each clause is a disjunction of literals, where each literal represents either a variable or its negation, ensuring that at least one of the included variables is assigned the required value for the clause to hold. All clauses in the \acrshort{cnf} must be satisfied in a valid assignment to the variables. This structure allows efficient constraint propagation and satisfiability checking \citep{biere2009handbook}. 

The arithmetic component of the \acrshort{smt} formulation consists of difference constraints, which enforce relationships between integer variables. For a given activity $a = (i,j)$ these constraints take the following form:

\begin{equation}\label{eq:diff_example}
    t_i - t_j \geq \delta_a \lor \neg q,
\end{equation}

Where \( t_i \) and \( t_j \) are event times, \( \delta_a \) is a minimum required separation, and \( q \) is a Boolean variable controlling whether the constraint must hold. This encoding allows conditional constraints, meaning that if \( q \) is set to \( \texttt{False} \), the constraint is effectively deactivated. In contrast, the Boolean term \( \lor \neg q \) is omitted for mandatory constraints that must always be respected.
For further background on the extension of \acrshort{sat} with difference constraints, we refer the reader to \citet{armando2004sat}.

Since the \acrshort{smt} encoding requires that all variables involved in Boolean operations are binary, we preprocess the input \acrshort{ean} to ensure that all activity durations \( \delta_a \) fall within the range \( [0, T] \), thus converting all periodicity variables \( k_a \) into binary representations. To achieve compliance with this condition, we employ the method proposed by \citet{Peeters2003CyclicOptimization}. This transformation splits any activity with a non-binary \( k_a \) into two or more activities, ensuring that each resulting activity can be represented with a binary \( k_a \).

\subsubsection{Boolean Encoding of the \acrfull{tfn}} 

Next, we describe how to encode the problem in \acrshort{smt}. Initially, we define some helper functions before delving into encoding the \acrshort{mip} as outlined in \cref{pesp:formulation}.

\begin{subequations}
\begin{align}
    & \texttt{encode-at-most-one}(\mathcal{Q}) := \bigwedge_{\forall q_i,q_j \in \mathcal{Q},~i < j } (\neg q_{i} \lor \neg q_{j}) \label{eq:helper_most_1} \\
    & \texttt{encode-at-least-one}(\mathcal{Q}) := \bigvee_{q \in \mathcal{Q}} q \label{eq:helper_least_1} \\
    & \begin{aligned} \texttt{encode-exactly-one}(\mathcal{Q}) := & \ ~~~\texttt{encode-at-least-one}(\mathcal{Q}) \\
    & \land \texttt{encode-at-most-one}(\mathcal{Q}) \end{aligned} \label{eq:helper_exact_1}
\end{align}
\end{subequations}

The \texttt{encode-at-most-one} function (\cref{eq:helper_most_1}) ensures that no more than one variable in a set is assigned the value \( \texttt{True} \). Conversely, \texttt{encode-at-least-one} (\cref{eq:helper_least_1}) ensures at least one variable in a set is assigned the value $\texttt{True}$. Combining these two functions, \texttt{encode-exactly-one} (\cref{eq:helper_exact_1}) guarantees that exactly one variable in a set is assigned the value $\texttt{True}$.

Next, we focus on encoding the \acrshort{tfn} and vehicle circulation cardinality constraints. This encoding step involves defining a \acrshort{sat} variable \(q \in \mathcal{Q}\) for each \(v \in \mathcal{V}\) node, \(w \in \mathcal{W}\) link, and $u \in \mathcal{U}$ vehicle circulation link in the \acrshort{tfn}.

\begin{subequations}
\begin{align}
\texttt{encode-TFN}~(\mathcal{V}) & := \bigwedge_{\forall v \in \mathcal{V}}(
\neg q_{v} \bigvee_{y \in \alpha^+(v)} (q_{y} 
\land \neg q_{v}) \bigvee_{y \in \alpha^-(v)} q_{y} \nonumber \\
& \land \texttt{encode-at-most-one}(\mathcal{Q}_{(\alpha^+(v))})  \nonumber\\
& \land \texttt{encode-at-most-one}(\mathcal{Q}_{(\alpha^-(v))})) \label{eq:encode-TFN} \\
\texttt{encode-one-train-path} & := \bigwedge_{\forall c \in \mathcal{C}} 
\texttt{encode-exactly-one}({q_v : v \in \mathcal{V}_c}) \label{eq:encode-one-train-path}
\end{align}
\end{subequations}

To encode the network, we can use the same approach as given by \citet{FuchsTrivellaCorman2022}, which encodes flow balance in \cref{eq:encode-TFN} and then requires one path per train in \cref{eq:encode-one-train-path}. 

\subsubsection{Encoding of the \acrfull{ean}}  

Having encoded the \acrshort{tfn}, we now turn our attention to encoding the \acrshort{ean} for the \acrshort{smt} equivalent of model~\eqref{pesp:formulation}. We begin by defining \acrshort{sat} variables \(\mathcal{Q}\) and time variables. Each event \(e \in \mathcal{E}\) is associated with a positive time variable \(t_e\). 

\begin{subequations}
\begin{align}
\texttt{encode-event}(\mathcal{E}) & := \bigwedge_{\forall e \in \mathcal{E}}  (t_e - t_0 \geq 0) 
\land (t_0 - t_e \geq -T) \label{eq:encode-event}
\end{align}
\end{subequations}

Each activity \( a = (i, j) \in \mathcal{A} \) links two events and has an associated period offset \( k_a \in \{0,1\} \), due to the preprocessing step described in the \acrshort{mip} formulation in Section~\ref{ss:mip_form}. This binary nature allows us to represent each activity using exactly two precedence directions: either \( t_j \geq t_i \), corresponding to \( k_a = 0 \), or \( t_j < t_i \), which implies \( k_a = 1 \). To encode this behavior, we introduce two Boolean variables for each activity: \( q_a \) to represent the case where \( t_j \geq t_i \), and \( \hat{q}_a \) to describe the inverted precedence \( t_i > t_j \). These two cases are mutually exclusive and collectively exhaustive, ensuring that precisely one is active when the activity is selected. The corresponding difference constraints are then conditionally enforced using the variables \( q_a \) and \( \hat{q}_a \), as shown in \cref{eq:encode-minimal-delta,eq:encode-maximal-delta}. Finally, \cref{eq:activate-precedence} ensures that precisely one of the two precedence directions is activated whenever the activity is selected (i.e., when all elements in \( \Omega_a \) are active).

\begin{subequations}
\begin{align}
\texttt{encode-minimal-delta}(\mathcal{A}) & := \bigwedge_{\forall a \in \mathcal{A}} 
((t_j - t_i \geq \delta^{\min}_a) \lor \neg q_a) 
\land ((t_i - t_j \leq (T - \delta^{\min}_a)) \lor \neg \hat{q_a}) \label{eq:encode-minimal-delta} \\
\texttt{encode-maximal-delta}(\mathcal{A}) & := \bigwedge_{\forall a \in \mathcal{A}} 
((t_i - t_j \geq -\delta^{\max}_a) \lor \neg q_a) 
\land ((t_j - t_i \leq -(T - \delta^{\max}_a)) \lor \neg \hat{q_a}) \label{eq:encode-maximal-delta} \\
\texttt{activate-precedence}(\mathcal{A}) & := \bigwedge_{\forall a \in \mathcal{A}} 
((q_a \lor \hat{q_a}) \lor \bigwedge_{\omega \in \Omega_a} \neg q_{\omega}) \label{eq:activate-precedence}
\end{align}
\end{subequations}

\subsubsection{Encoding the Vehicle Circulation}

After encoding the \acrshort{tfn} and \acrshort{ean}, we need to encode the vehicle circulation before focusing on transforming the encoded satisfiability problem into a minimization task. Therefore, we enforce that each start event \(e \in \mathcal{E}_{\text{Start}}\) and each end event \(e \in \mathcal{E}_{\text{End}}\) is connected by exactly one transfer arc:

\begin{subequations}
\begin{align}
\texttt{encode-vehicle-circulation-start}(\mathcal{E}_{\text{Start}}) 
& := \bigwedge_{\forall e \in \mathcal{E}_{\text{Start}}} 
\texttt{encode-exactly-one}(\mathcal{Q}_{(\alpha^+(e))}) \label{eq:encode-vehicle-circulation-start} \\
\texttt{encode-vehicle-circulation-end}(\mathcal{E}_{\text{End}}) 
& := \bigwedge_{\forall e \in \mathcal{E}_{\text{End}}} 
\texttt{encode-exactly-one}(\mathcal{Q}_{(\alpha^-(e))}) \label{eq:encode-vehicle-circulation-end}
\end{align}
\end{subequations}

These \cref{eq:encode-vehicle-circulation-start,eq:encode-vehicle-circulation-end} ensure the correct transfers of vehicles by enforcing flow conservation and exclusivity at transfer events, thereby guaranteeing that each scheduled train service is connected to a feasible vehicle cycle, consistent with the operational rules defined in the \acrshort{smt} model.

\begin{equation}
    \texttt{encode-vehicle-count}(n, \mathcal{S} \cup \mathcal{U}) := 
    \texttt{sequence-counter}(n, \hat{P})
    \quad \text{where} \quad \hat{P} = \{\hat{p}_{a} \mid a \in \mathcal{S} \cup \mathcal{U}\} \label{eq:encode-vehicle-count}
\end{equation}

As outlined in Section~\ref{ss:mip_form}, the number of vehicles required to operate a feasible periodic timetable corresponds to the sum over all \( k_a \in \{0, 1\} \) for selected commercial and vehicle-transfer activities \( a \in \mathcal{S} \cup \mathcal{U} \). In the \acrshort{smt} formulation, each binary variable \( \hat{p}_a \) encodes whether the corresponding precedence is inverted, i.e., whether \( k_a = 1 \). Consequently, \cref{eq:encode-vehicle-count} applies the \texttt{sequence-counter} encoding \citep{SeqCounterCarstenSinz2005} to the set of literals \( \hat{p}_a \), thereby enforcing an upper bound \( n \) on the number of vehicles. This reformulates the vehicle circulation problem as a feasibility check under a fixed fleet size constraint.

\subsection{Translation to SAT}\label{sec:sat_translation}

To enable a direct comparison with the \acrshort{smt} and \acrshort{mip} models, we formulate a \acrshort{sat} encoding of the periodic timetabling problem. This encoding follows the approach of \citet{FuchsTrivellaCorman2022}, which extends the method of \citet{Gromann2011PolynomialSATb} to handle train routing, and we further adapt it to incorporate vehicle circulation constraints. 

For brevity, we do not explicitly detail the encoding of train sequencing, headway constraints, or vehicle circulation, as these follow the structure already established in \citet{FuchsTrivellaCorman2022}. Instead, we outline the key aspects distinguishing the \acrshort{sat} model from \acrshort{smt} and \acrshort{mip}. 

The transformation consists of two primary steps. First, all event times are encoded using an \textit{order encoding}, representing integer time values as Boolean variables. Second, all difference constraints—previously formulated in \acrshort{smt}—are expressed in propositional logic. The resulting model ensures consistency across train movements while maintaining routing flexibility. 

Vehicle circulation constraints are incorporated analogously to the \acrshort{smt} model by enforcing \texttt{exactly-one} constraints for train handovers at terminal stations: each terminating train selects precisely one outgoing vehicle-transfer activity, and each originating train selects exactly one incoming activity. To encode the vehicle count, we introduce indicator literals \( \hat{p}_a \) for each vehicle-related activity \( a \in \mathcal{S} \cup \mathcal{U} \), representing whether the activity uses the cyclic wrap-around (i.e., \( k_a = 1 \)). These literals are used in conjunction with a \texttt{sequence-counter} encoding \citep{SeqCounterCarstenSinz2005} to enforce an upper bound \( n \) on the number of wrap-arounds—and hence vehicles. This reformulates the vehicle circulation problem as a feasibility check under a fixed fleet size constraint, consistent with the \acrshort{smt} model.

\section{Results}\label{sec:results}

The results presented in this section are based on three sets of experiments designed to evaluate the efficiency of the proposed modeling approaches and assess the impact of various problem characteristics. First, we describe the implementation of the models (Section~\ref{sec:implementation}), and present the case-study and the set of instances derived from the Swiss railway network (Section~\ref{ss:case_infra}), and compare the computational performance of Mixed Integer Programming (\acrshort{mip}), Boolean Satisfiability (\acrshort{sat}), and Satisfiability Modulo Theories (\acrshort{smt}) solvers varying the time discretisation granularity (Section~\ref{ss:performance_comparison}). Second, we investigate how discretisation affects vehicle requirements under two circulation strategies: vehicles restricted to the same line versus shared use across all trains (Section~\ref{ss:effect_discretisation}). Third, we analyze the effect of train routing flexibility on vehicle counts, comparing fixed and flexible routing scenarios with the same two circulation strategies (Section~\ref{ss:impact_fixed_routes}).
Finally, we study the benefit of the proposed integrated problem against sequential equivalent planning procedures (Section~\ref{ss:benefit_of_integration}).

\subsection{Implementation}\label{sec:implementation}

The proposed models are implemented using the \textit{OpenBus} framework \citep{an_open_toolbox}, ensuring consistency across all solving methodologies. We consider three optimization approaches: \acrshort{mip}, \acrshort{sat}, and \acrshort{smt}, each using a computing server equipped with four CPU cores (Intel Xeon Gold 6248) and 32 GB of RAM for all computational experiments. To leverage parallelization, all solvers utilize four threads.

The \acrshort{mip} formulation is implemented using Gurobi 12.0.1 \citep{gurobi}, with four solver threads assigned. The \acrshort{sat} formulation is implemented using Glucose 4.1 \citep{glucose_4_1} via the PySAT package \citep{imms-sat18}. A portfolio strategy utilizes all four cores, where each core runs an independent solver instance initialized with a different random seed \citep{horde_sat}. The \acrshort{smt} solver extends the approach of \citet{LEUTWILER2022525}, employing a portfolio-based strategy similar to the \acrshort{sat} approach. The solvers do not share any state or information, as they work on independent search spaces, with the first to terminate providing the final result.

We solve the \acrshort{sat} and \acrshort{smt} problems using an \textit{ascending linear search} to determine the minimal vehicle count. First, a lower bound on the number of required vehicles is computed based on relaxed circulation constraints, neglecting headway constraints. This relaxation provides an initial lower bound for subsequent iterations. The model is then solved incrementally, starting from this bound and increasing the vehicle count \( n \) step by step. If the model is infeasible for a given \( n \), the vehicle count is incremented by one, and the lower bound is updated until feasibility is attained. Once a feasible solution is found, it is guaranteed optimal, as during this linear search, all instances with fewer vehicle counts have been proven to be infeasible. Thus, we can conclude the procedure.

\subsection{Infrastructure and Instances}\label{ss:case_infra}

For our case study, we used data provided by \acrfull{rhb}, a Swiss railway company operating most of the railway lines of the canton of Grisons. Following a methodology similar to \citet{FuchsTrivellaCorman2022}, we modeled the network at a mesoscopic level, as many sections consist of a single track. The infrastructure spans 380~km, with a complex terrain and operational restrictions. Technical running times were calculated using the same procedures as those employed by \acrshort{rhb}, ensuring that the instances reflect realistic railway operations.

\begin{figure}[H]%
\begin{center}
\includegraphics[width=0.8\textwidth]{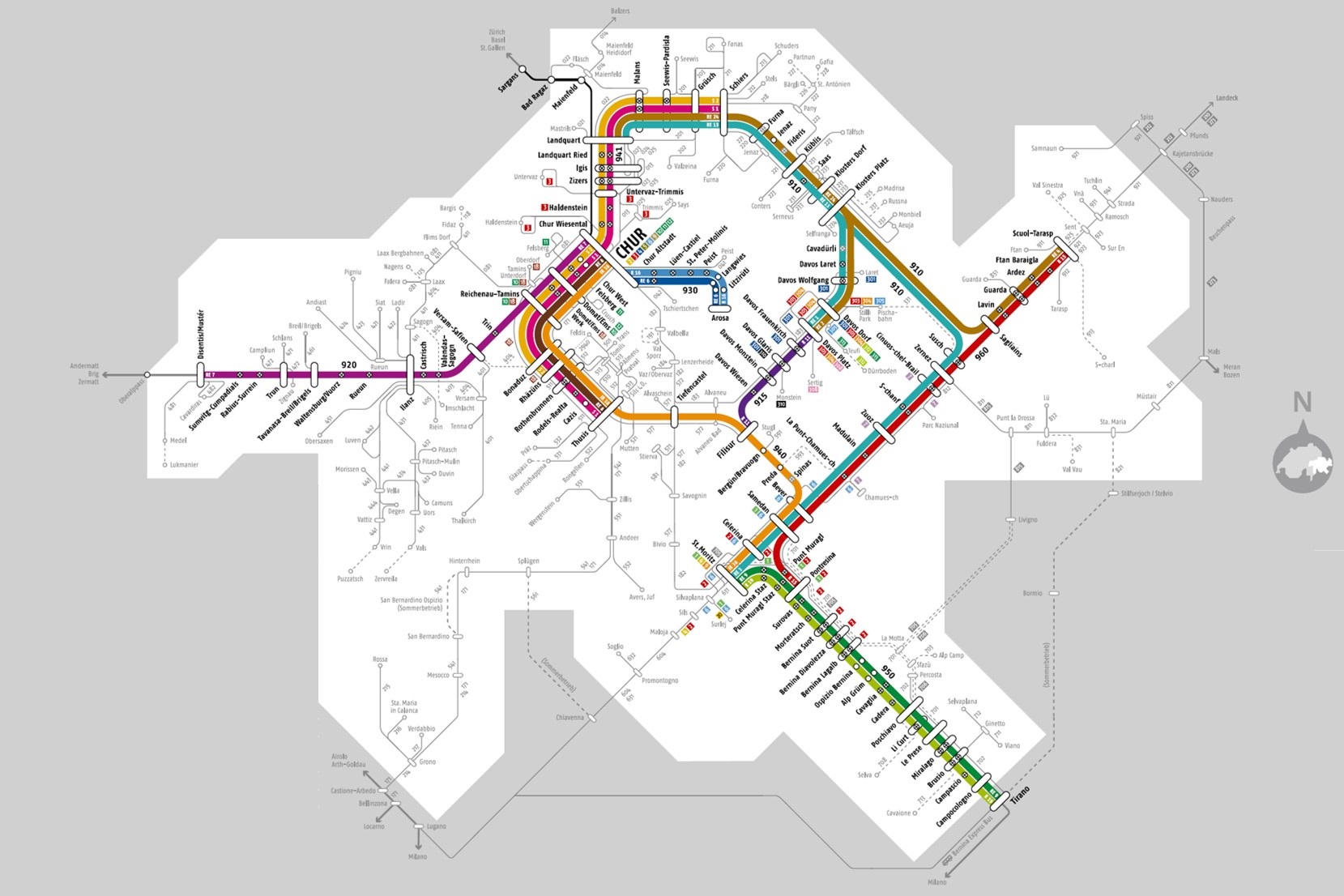}%
\caption{The current line plan for \acrshort{rhb} \citep{Streckennetz_rhb}.}%
\label{fig:infrastructure_in_detail}%
\end{center}
\end{figure}

Using the current line plan (see \cref{fig:infrastructure_in_detail}), we generated a series of problem instances ranging from 1 to 10 lines, ensuring that the selected line plans remained connected. The characteristics of these instances, including the number of events, activities, headway constraints, and routing alternatives, are summarized in \cref{tab:instance_characteristics}. 

Each instance represents an increasing level of complexity, with additional lines introducing more routing alternatives, activities, and constraints. The number of vehicle transfer links varies depending on the vehicle-sharing policy:
\begin{itemize}
    \item \texttt{No-Sharing}: Vehicles remain restricted to operating within their assigned service, meaning they cannot transfer between different lines.
    \item \texttt{Full-Sharing}: Vehicles can be shared across different lines (without dead-heading), allowing for a more flexible assignment and potentially reducing the number of required vehicles.
\end{itemize}

\begin{table}[H]
\centering
\caption{Instance characteristics by number of lines.}
\label{tab:instance_characteristics}
\resizebox{0.7\textwidth}{!}{%
\begin{tabular}{c|c|c|c|c|c|c|c}
\toprule
Lines & Events & Activities & Headway & Itinerary & Routing & \multicolumn{2}{c}{Vehicle Transfers} \\
      &        &            & Constraints & Activities & Alternatives & \texttt{~~No~~} & \texttt{Full} \\
\midrule
1  & 160  & 1038   & 812   & 144  & 34   & 2  & 2  \\
2  & 424  & 2799   & 2234  & 385  & 95   & 4  & 6  \\
3  & 608  & 4316   & 3515  & 541  & 124  & 6  & 12 \\
4  & 1011 & 12317  & 10986 & 889  & 196  & 8  & 14 \\
5  & 1363 & 22134  & 20302 & 1232 & 287  & 10 & 22 \\
6  & 1975 & 38640  & 36059 & 1753 & 405  & 12 & 32 \\
7  & 2379 & 43343  & 40243 & 2114 & 490  & 14 & 36 \\
8  & 2891 & 53404  & 49681 & 2535 & 573  & 16 & 42 \\
9  & 2989 & 54339  & 50478 & 2619 & 590  & 18 & 46 \\
10 & 3291 & 58425  & 54144 & 2887 & 644  & 20 & 50 \\
\bottomrule
\end{tabular}%
}
\end{table}

As shown in \cref{tab:instance_characteristics}, the number of events and activities increases with the number of lines, naturally resulting in more headway constraints, itinerary activities (dwells and trips), and routing alternatives. Vehicle transfer options increase with instance size, especially when allowing vehicles to circulate across different lines. This flexibility, denoted as \textit{Full} in the table, is expected to reduce the overall fleet size compared to the \textit{No} case, where vehicles are restricted to individual lines.

\subsection{Comparison of Performance}\label{ss:performance_comparison}

To evaluate the performance of the three approaches—Mixed Integer Programming (\acrshort{mip}), Boolean Satisfiability (\acrshort{sat}), and Satisfiability Modulo Theories (\acrshort{smt})—we solve test instances derived from subsets of the \acrshort{rhb} line plan. Each instance in \cref{tab:instance_characteristics} is computed once under the \texttt{Full-Sharing} and \texttt{No-Sharing} policies. To assess the impact of time granularity, we solve each instance at four different discretisation levels: 6, 3, 2, and 1 seconds. To account for performance variability, each solver is executed five times per instance with different random seeds and a time limit of 5 hours. The plots below show the median computation times per approach. We also indicate the timeout threshold and highlight scaling behavior as the instance size increases. 

\begin{figure}[H]
    \centering
    \begin{subfigure}[b]{0.49\textwidth}
        \centering
        \includegraphics[width=\textwidth]{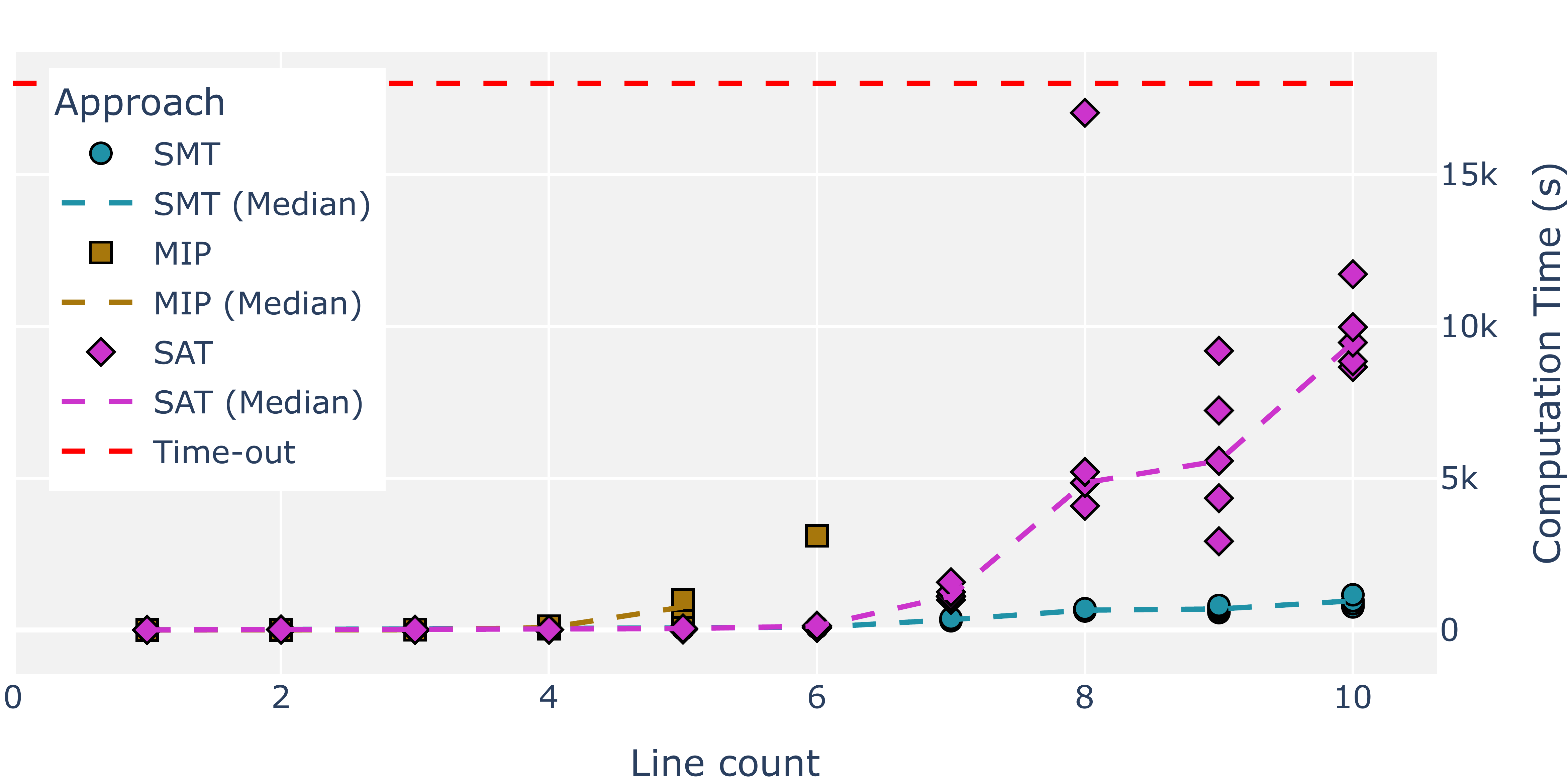}
        \caption{6-second granularity.}
        \label{fig:comp_time_dt_6_ns}
    \end{subfigure}
    \hfill
    \begin{subfigure}[b]{0.49\textwidth}
        \centering
        \includegraphics[width=\textwidth]{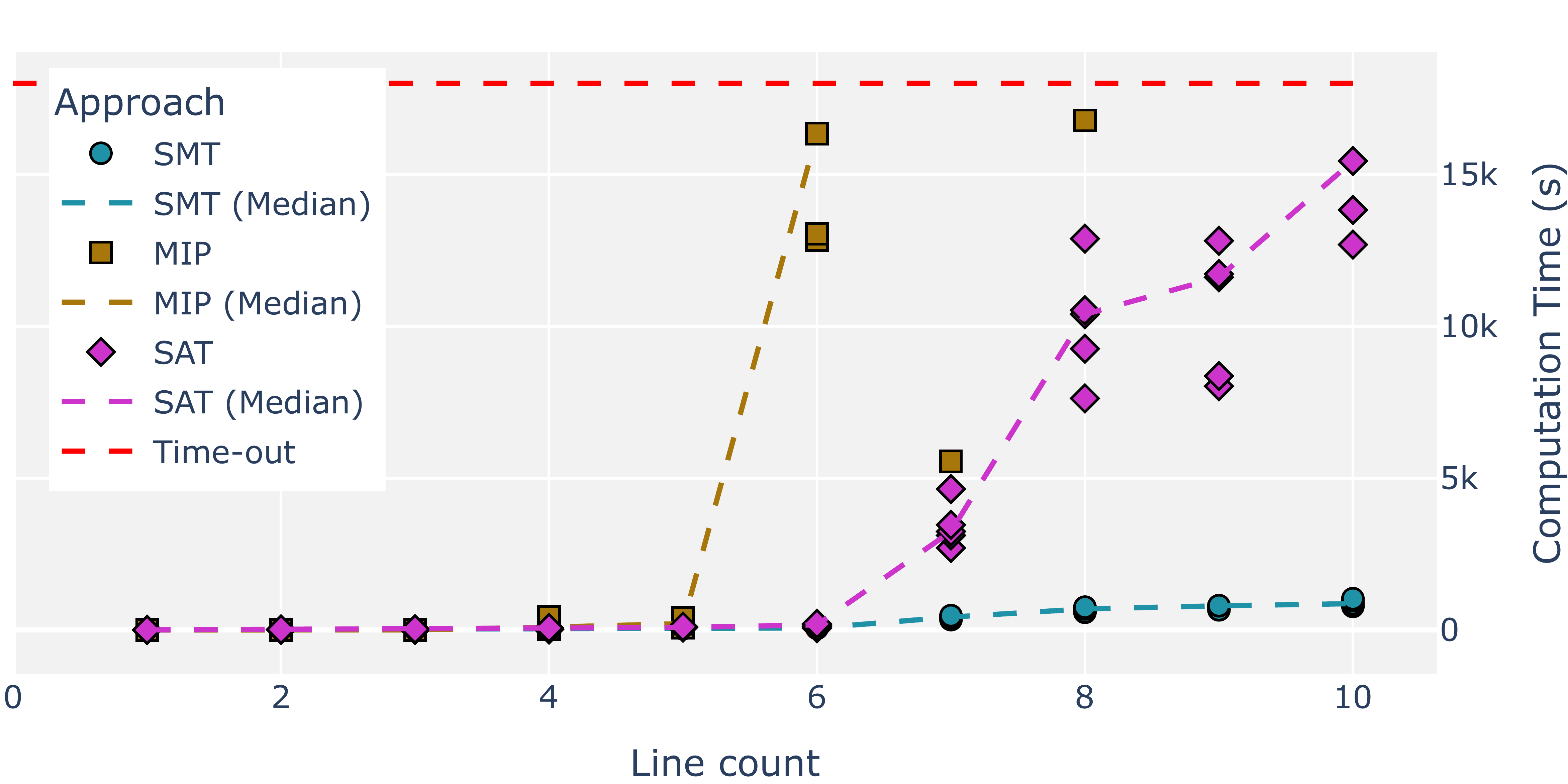}
        \caption{3-second granularity.}
        \label{fig:comp_time_dt_3_ns}
    \end{subfigure}
    \begin{subfigure}[b]{0.49\textwidth}
        \centering
        \includegraphics[width=\textwidth]{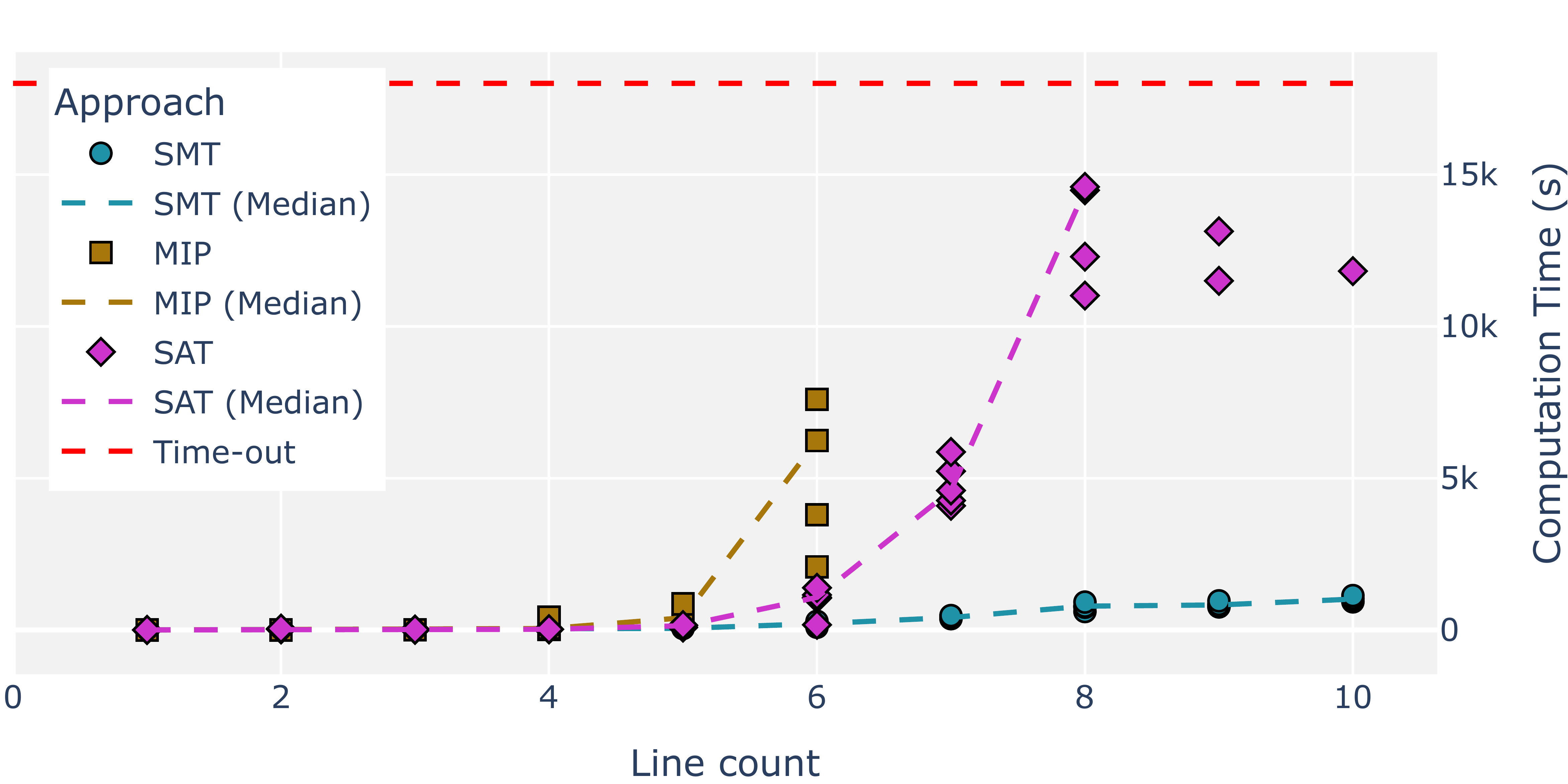}
        \caption{2-second granularity.}
        \label{fig:comp_time_dt_2_ns}
    \end{subfigure}
    \hfill
    \begin{subfigure}[b]{0.49\textwidth}
        \centering
        \includegraphics[width=\textwidth]{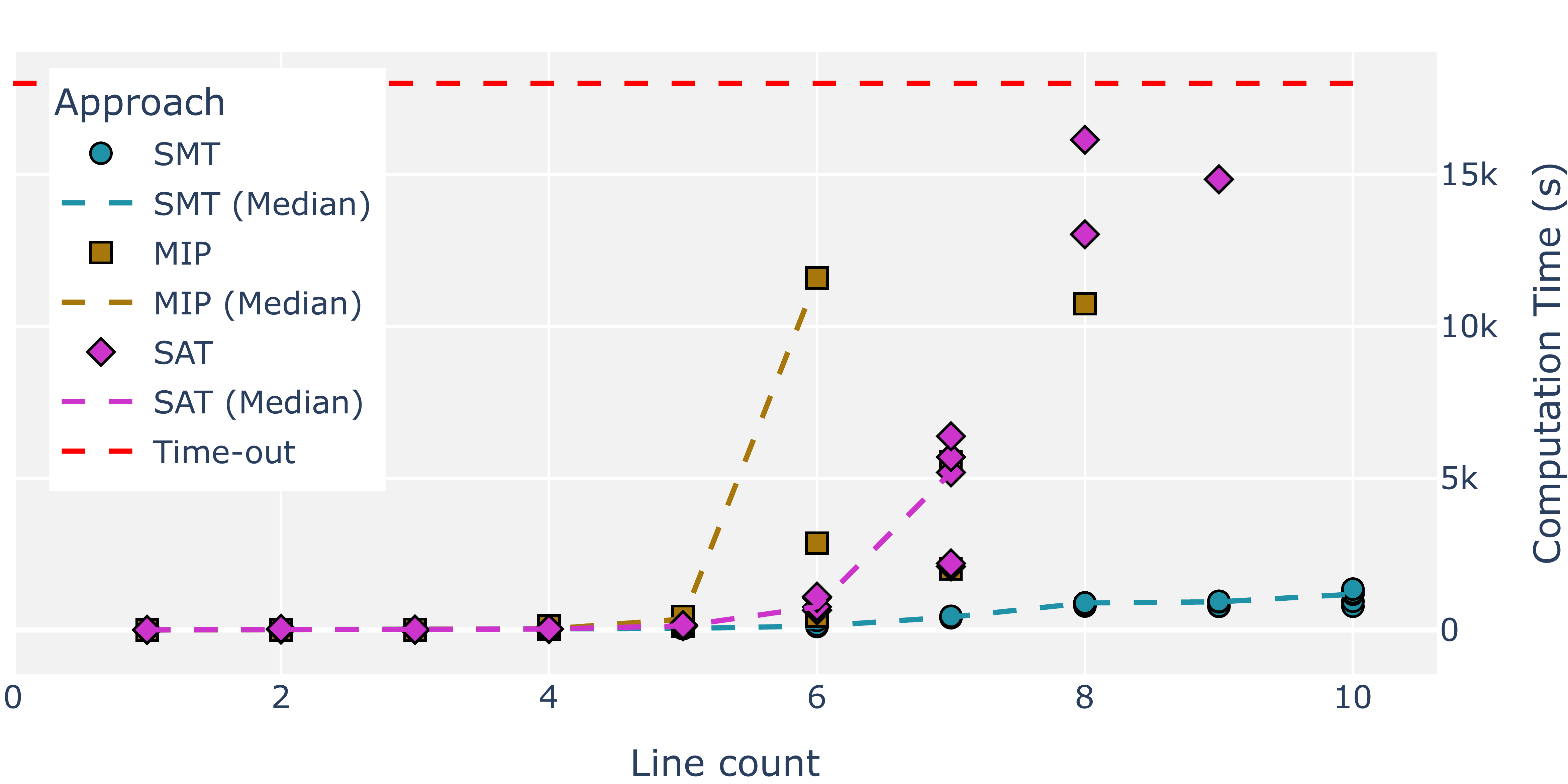}
        \caption{1-second granularity.}
        \label{fig:comp_time_dt_1_ns}
    \end{subfigure}
    \caption{Computation times for \acrshort{mip}, \acrshort{sat}, and \acrshort{smt} under \texttt{No-Sharing} strategy across four discretisation levels.}
    \label{fig:comp_time_no_sharing}
\end{figure}

The results under the \texttt{No-Sharing} policy in \cref{fig:comp_time_no_sharing} clearly demonstrate the superior scalability of the \acrshort{smt} formulation. While \acrshort{sat} and \acrshort{mip} exhibit acceptable performance at coarse resolutions (6 and 3 seconds), both degrade significantly as the temporal resolution increases. \acrshort{sat} fails to solve many instances beyond six lines at 1-second resolution, and \acrshort{mip} times out already at intermediate sizes. In contrast, \acrshort{smt} maintains stable runtime across all tested resolutions and solves all instances up to ten lines without reaching the time limit. This advantage in scalability becomes more pronounced as instance complexity grows, underlining the practical advantage of the \acrshort{smt} approach for large-scale, high-resolution periodic timetable optimization.

Under the \texttt{No-Sharing} policy in \cref{fig:comp_time_no_sharing}, we observe clear tipping points beyond which solvers fail to compute solutions within the time limit. At coarser time steps (6 and 3 seconds), \acrshort{sat} performs comparably well and outperforms \acrshort{mip}, solving all instances quickly. As the granularity increases, its performance degrades sharply. At the 1-second level, it frequently times out beyond six lines. \acrshort{mip} only handles networks up to 4–5 lines reliably across all resolutions and for these counts offers competitive runtime. However, for instances of larger sizes, \acrshort{mip} is no longer suitable. In contrast, \acrshort{smt} remains the most stable, showing consistent performance at high resolution and with larger networks. Notably, the difference in behavior between solvers is already visible at intermediate sizes, suggesting a gradual rather than sudden breakdown.

\begin{figure}[H]
    \centering
    \begin{subfigure}[b]{0.49\textwidth}
        \centering
        \includegraphics[width=\textwidth]{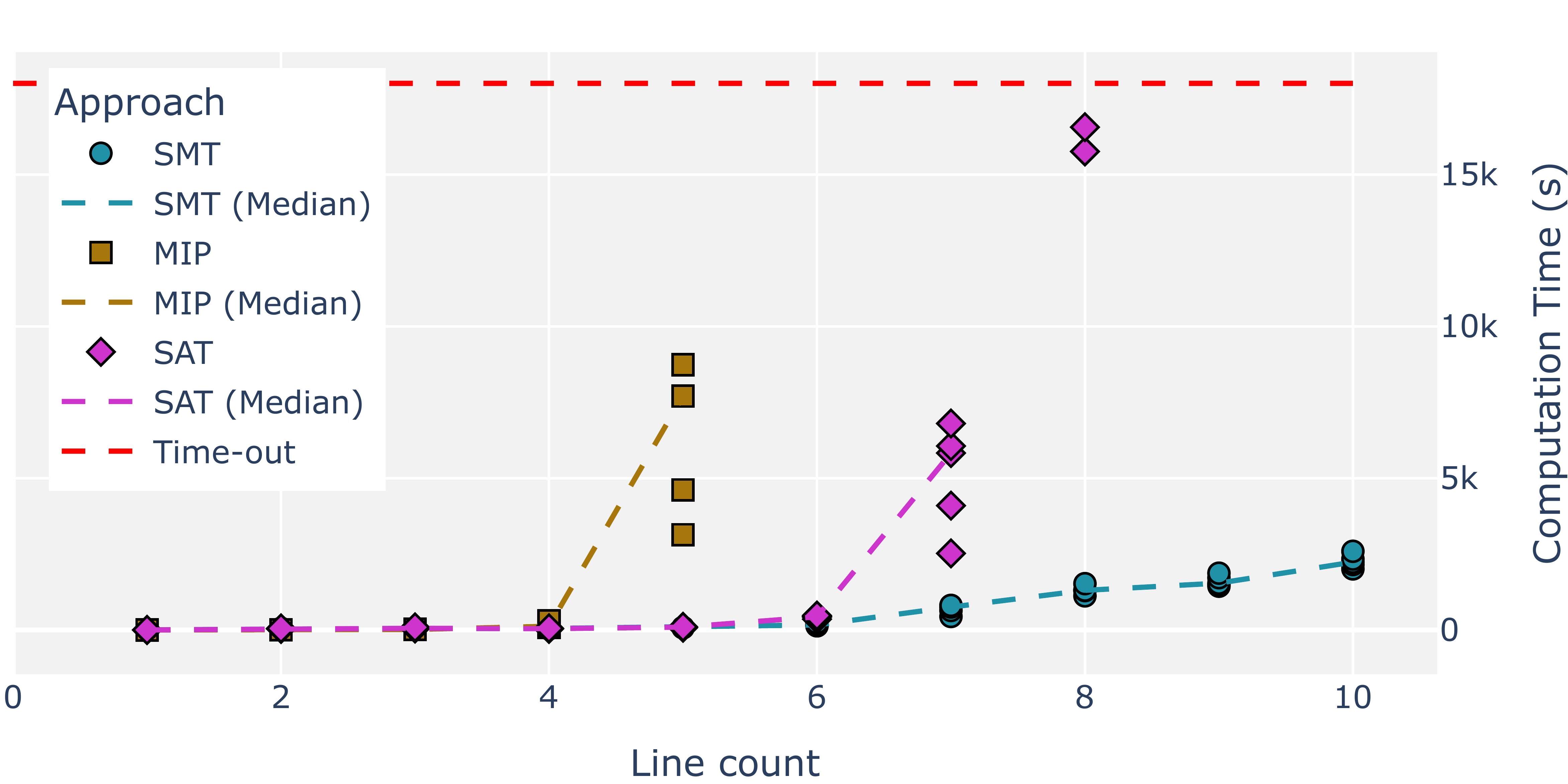}
        \caption{6-second granularity.}
        \label{fig:comp_time_dt_6_fs}
    \end{subfigure}
    \hfill
    \begin{subfigure}[b]{0.49\textwidth}
        \centering
        \includegraphics[width=\textwidth]{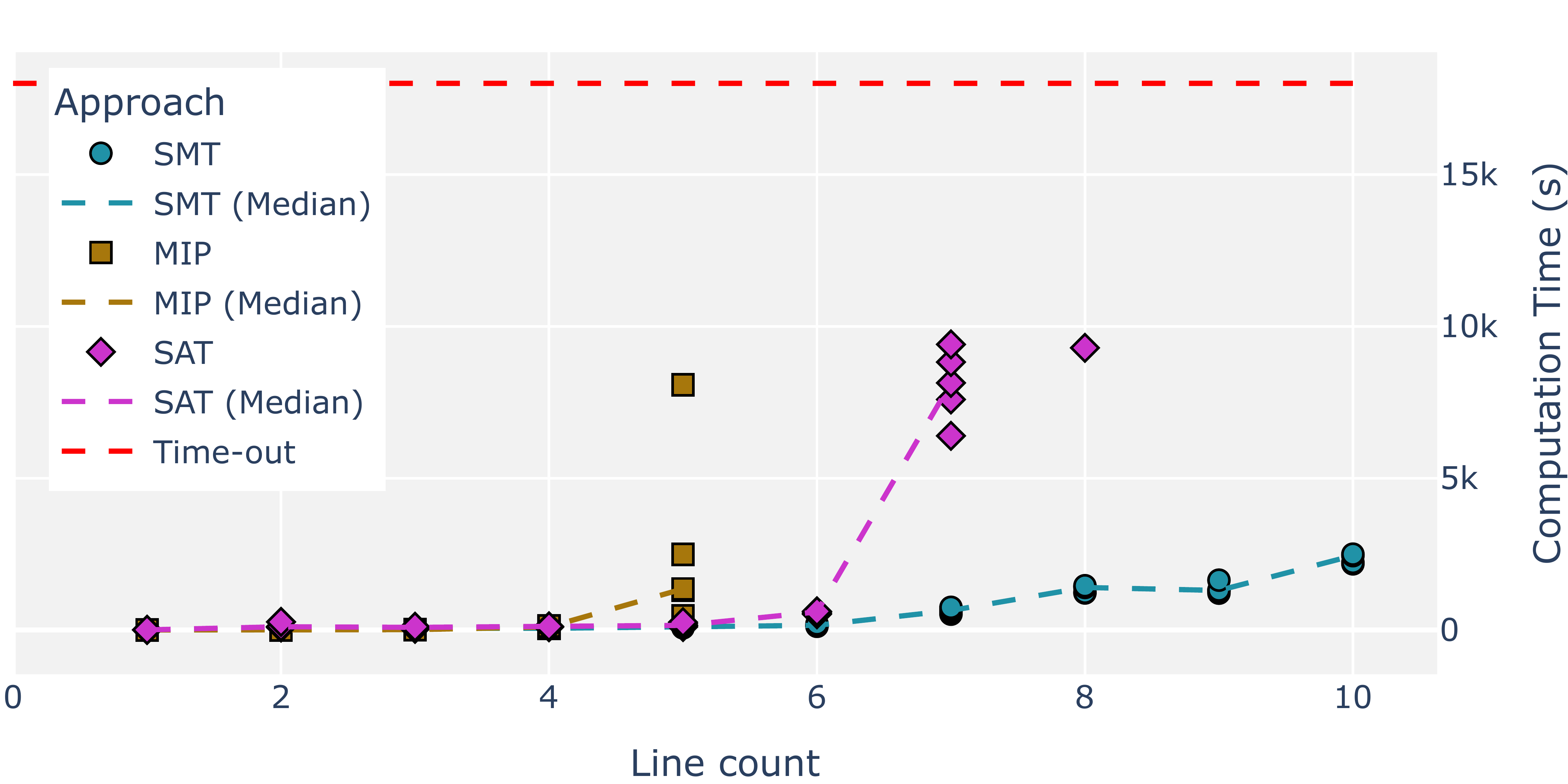}
        \caption{3-second granularity.}
        \label{fig:comp_time_dt_3_fs}
    \end{subfigure}
    \begin{subfigure}[b]{0.49\textwidth}
        \centering
        \includegraphics[width=\textwidth]{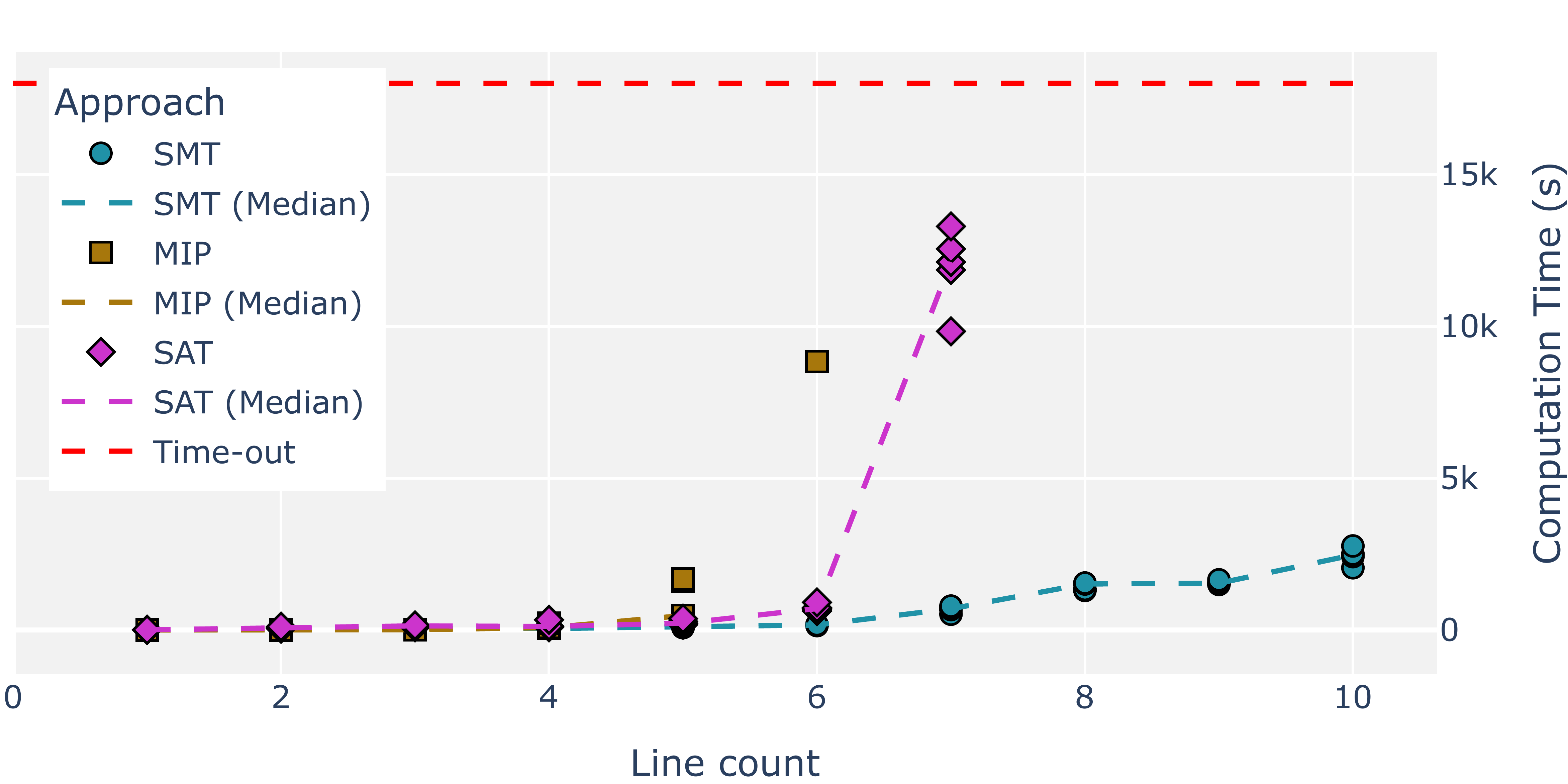}
        \caption{2-second granularity.}
        \label{fig:comp_time_dt_2_fs}
    \end{subfigure}
    \hfill
    \begin{subfigure}[b]{0.49\textwidth}
        \centering
        \includegraphics[width=\textwidth]{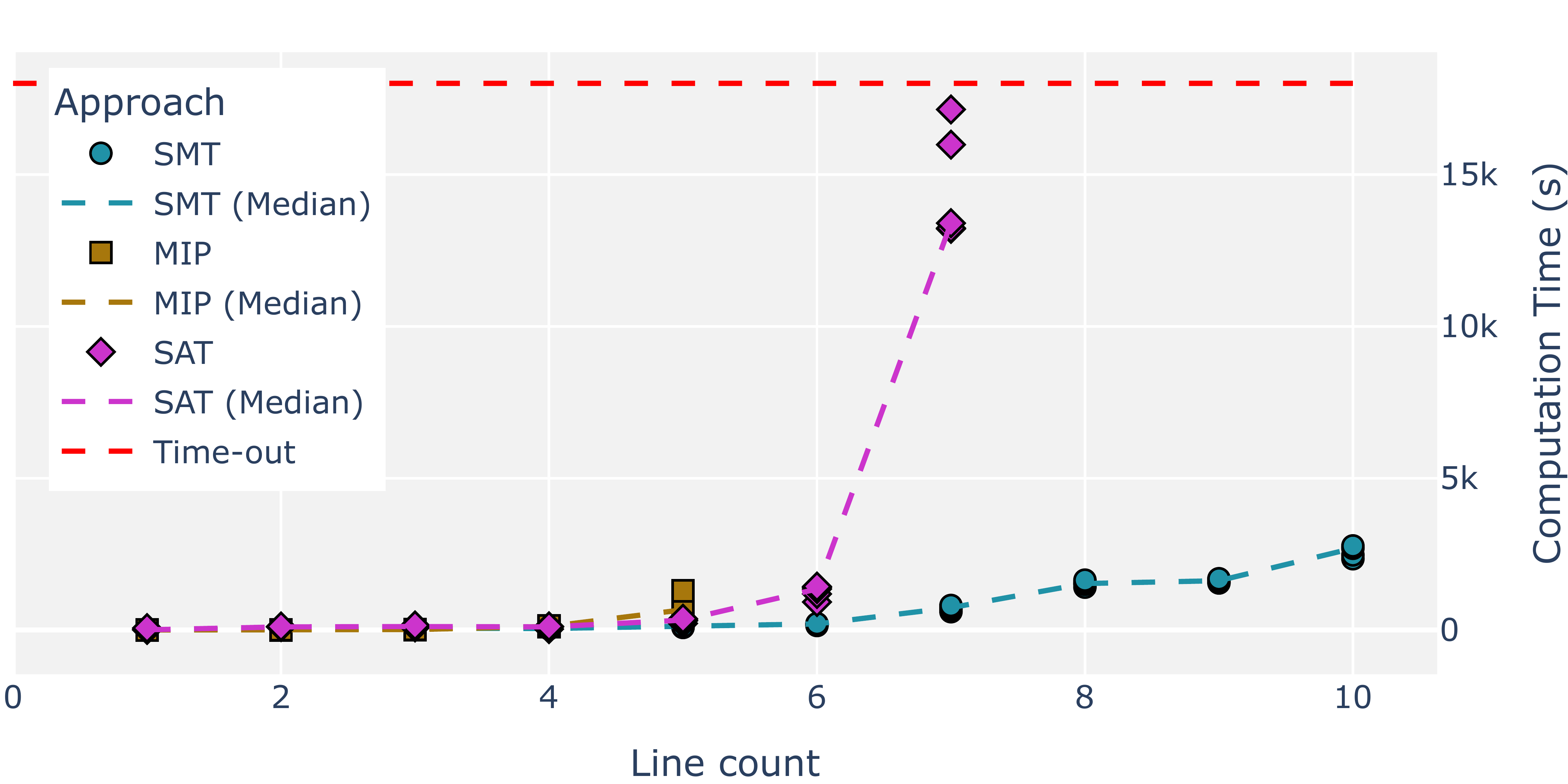}
        \caption{1-second granularity.}
        \label{fig:comp_time_dt_1_fs}
    \end{subfigure}
    \caption{Computation times for \acrshort{mip}, \acrshort{sat}, and \acrshort{smt} under \texttt{Full-Sharing} strategy across four discretisation levels.}
    \label{fig:comp_time_full_sharing}
\end{figure}

The patterns under \texttt{Full-Sharing} in \cref{fig:comp_time_full_sharing} are qualitatively similar but more pronounced than the ones with \texttt{No-Sharing} in \cref{fig:comp_time_no_sharing}. While \acrshort{sat} again performs well at coarse resolutions, it fails even earlier at finer ones. For example, at 1-second resolution, partial timeouts already appear from 5 lines onward. \acrshort{mip} behaves comparably to the \texttt{No-Sharing} case, while \acrshort{smt} continues to scale reliably. \texttt{Full-Sharing} appears to amplify the runtime demands of solvers, likely because solution space flexibility introduces additional combinatorial complexity.

To complement the runtime analysis, we report the cumulative success rate of each solver. We track whether a solver found a feasible solution for every instance and time limit and plot the fraction of solved instances over time. These plots provide a comprehensive overview of how quickly and reliably each solver performs across different instance sizes.

\begin{figure}[H]
    \centering
    \begin{subfigure}[b]{0.24\textwidth}
        \centering
        \includegraphics[width=\textwidth]{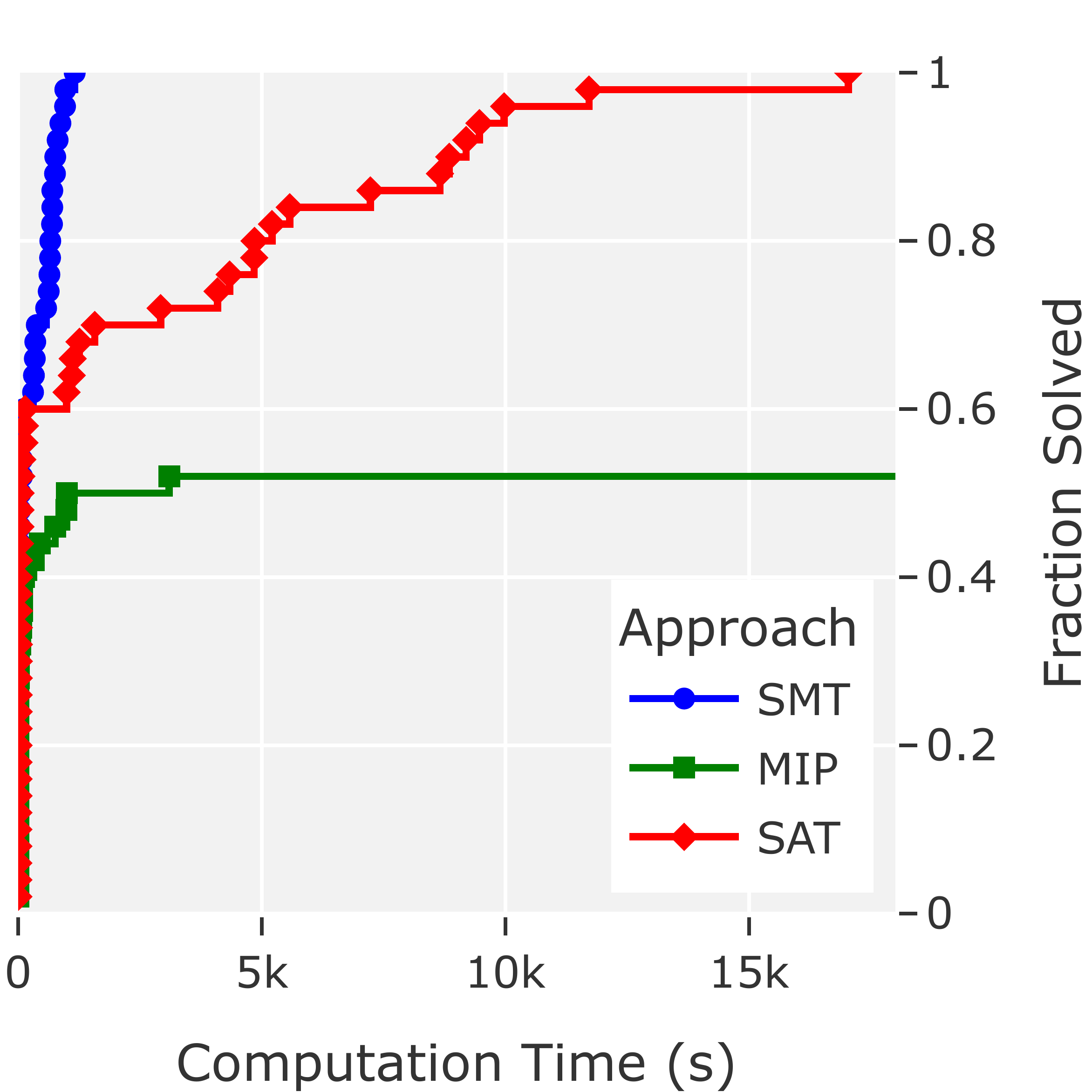}
        \caption{6-second granularity}
        \label{fig:cumulative_success_dt6_ns}
    \end{subfigure}
    \hfill
    \begin{subfigure}[b]{0.24\textwidth}
        \centering
        \includegraphics[width=\textwidth]{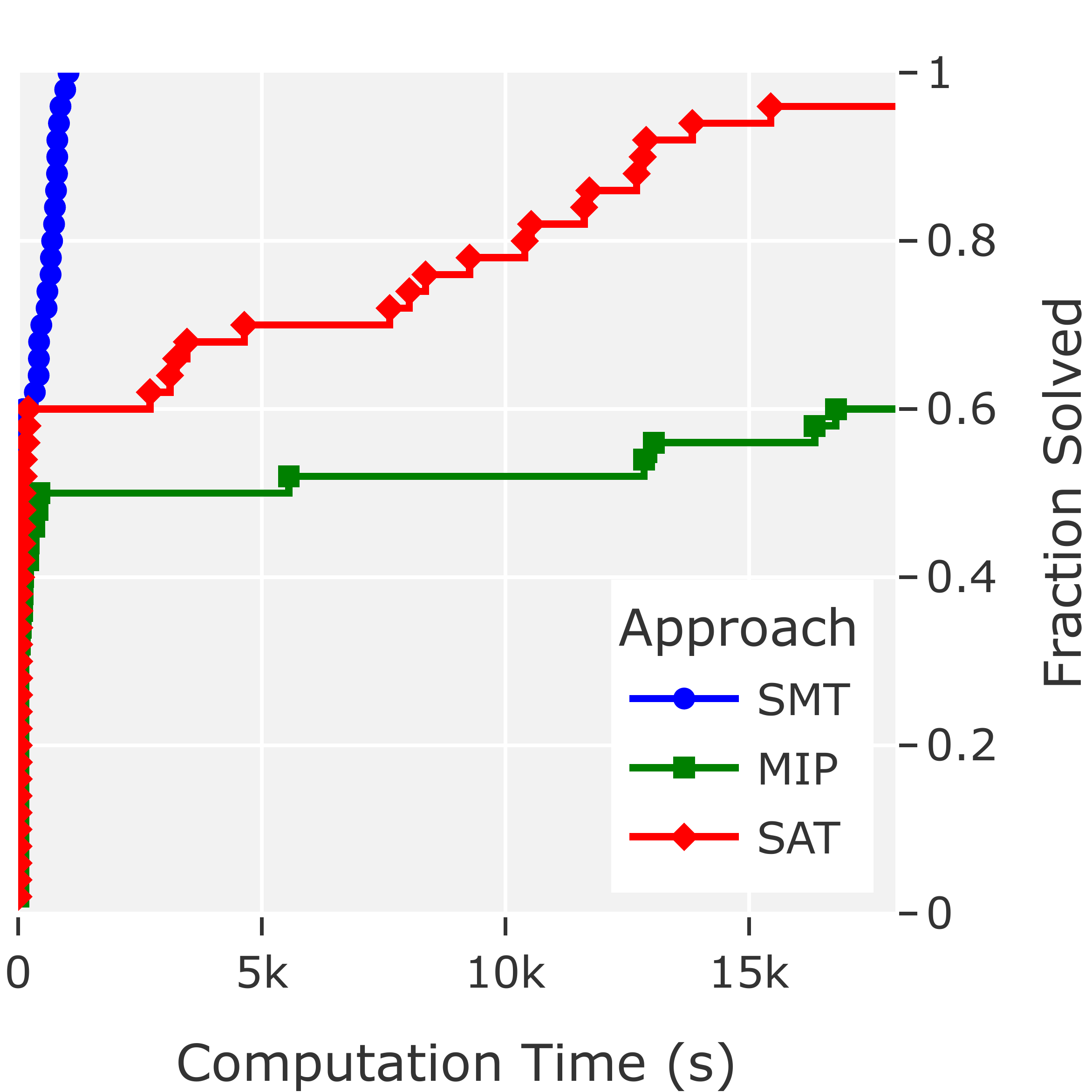}
        \caption{3-second granularity}
        \label{fig:cumulative_success_dt3_ns}
    \end{subfigure}
    \hfill
    \begin{subfigure}[b]{0.24\textwidth}
        \centering
        \includegraphics[width=\textwidth]{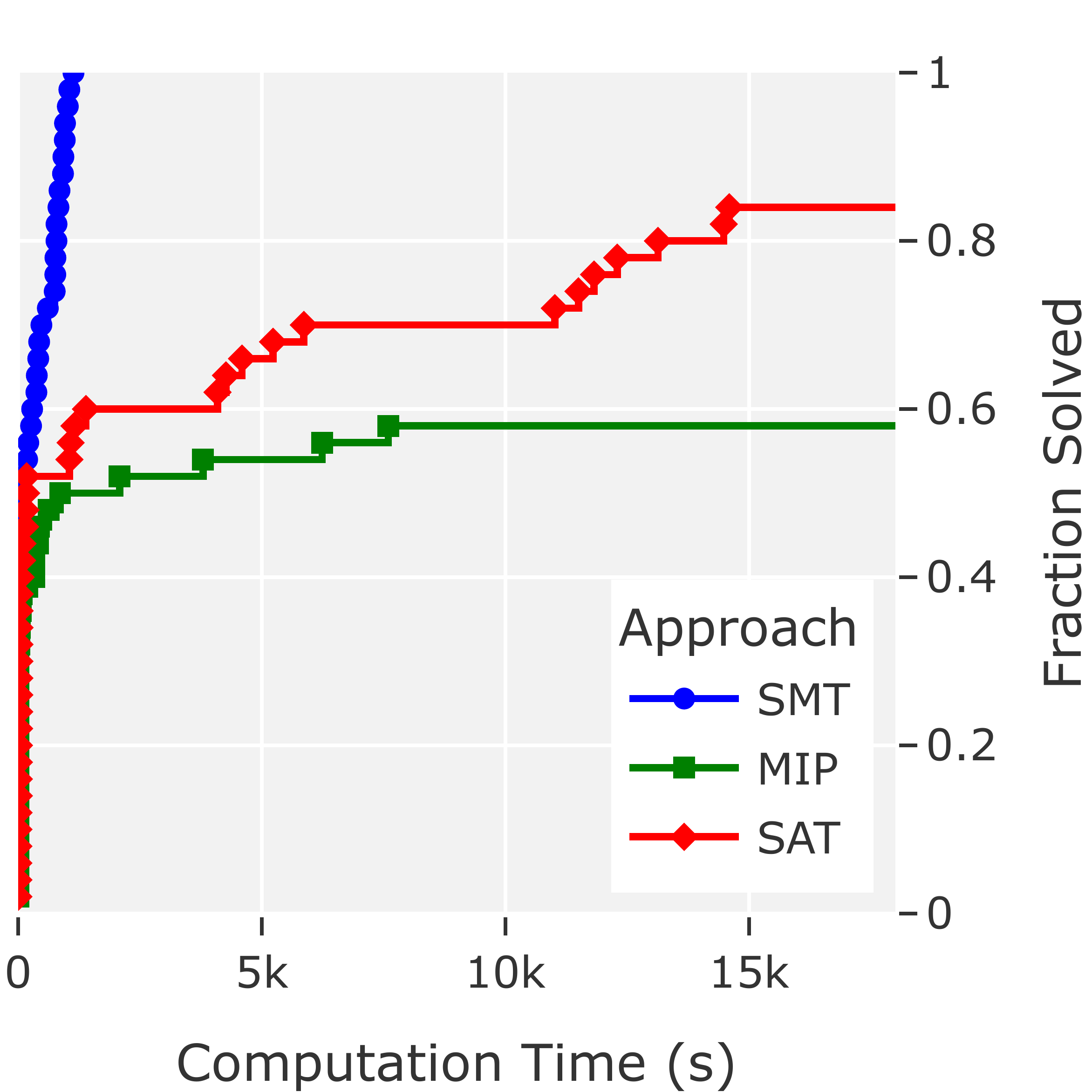}
        \caption{2-second granularity}
        \label{fig:cumulative_success_dt2_ns}
    \end{subfigure}
    \hfill
    \begin{subfigure}[b]{0.24\textwidth}
        \centering
        \includegraphics[width=\textwidth]{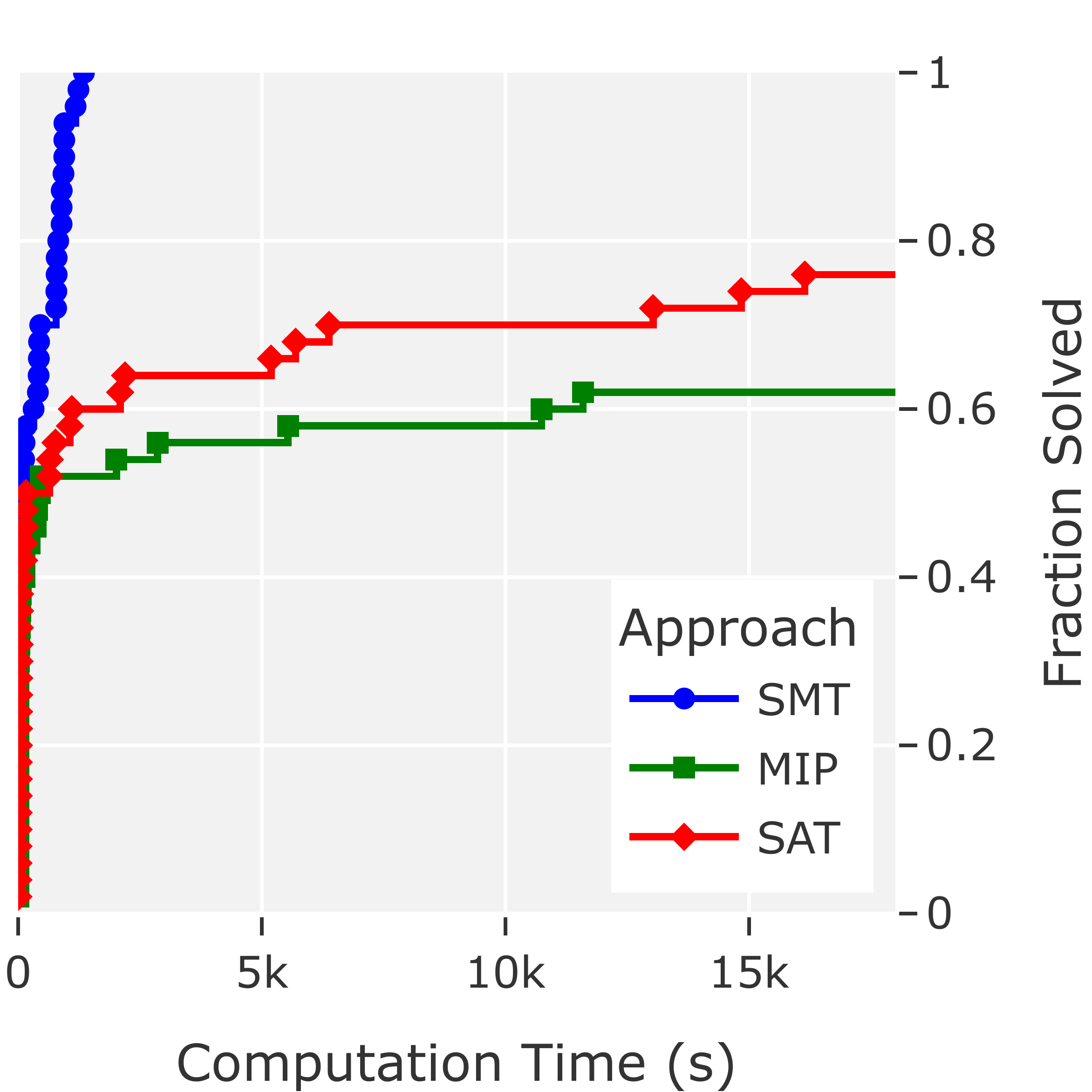}
        \caption{1-second granularity}
        \label{fig:cumulative_success_dt1_ns}
    \end{subfigure}
    \caption{Cumulative success rate for \acrshort{mip}, \acrshort{sat}, and \acrshort{smt} under \texttt{No-Sharing} strategy.}
    \label{fig:cumulative_success_ns}
\end{figure}

\begin{figure}[H]
    \centering
    \begin{subfigure}[b]{0.24\textwidth}
        \centering
        \includegraphics[width=\textwidth]{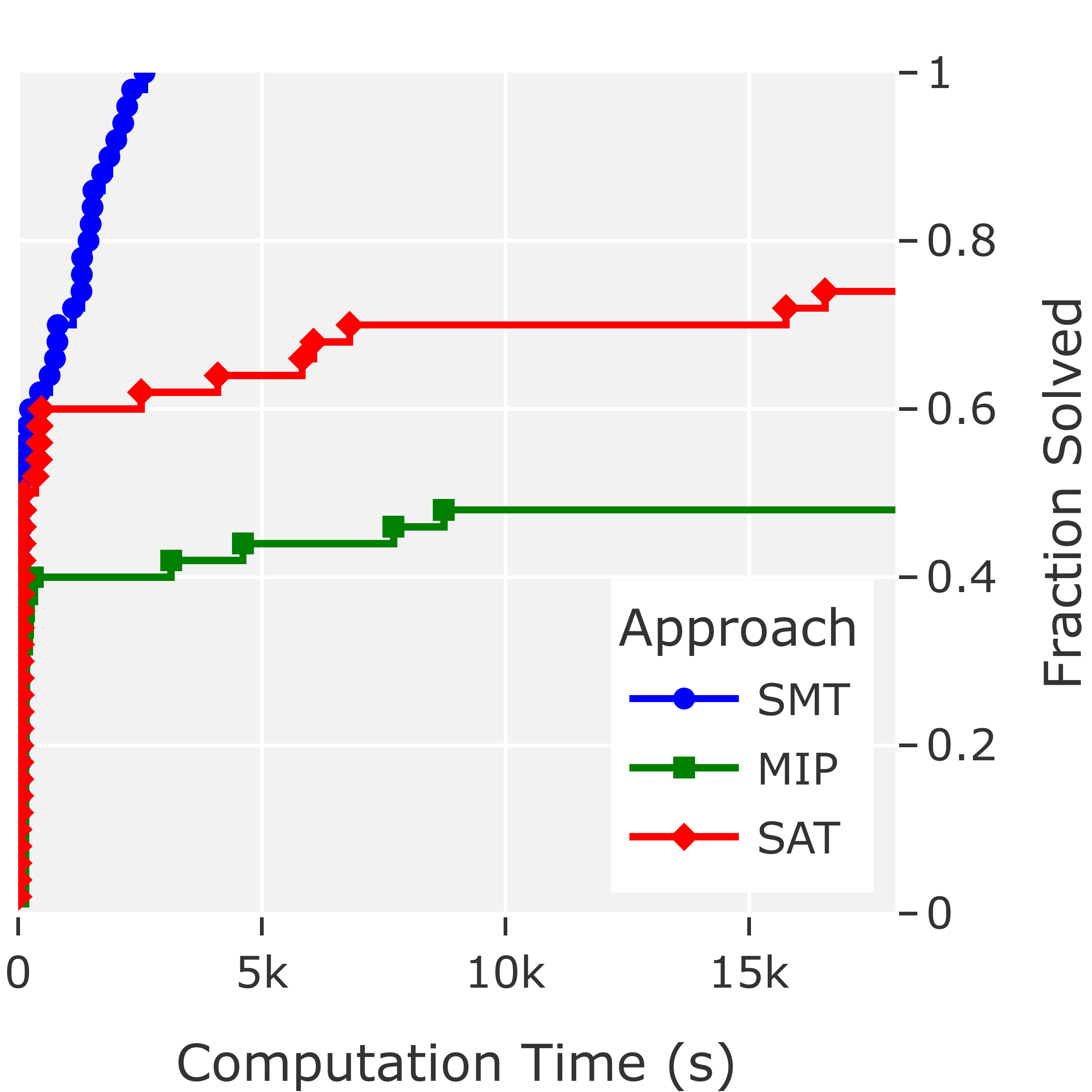}
        \caption{6-second granularity}
        \label{fig:cumulative_success_dt6_fs}
    \end{subfigure}
    \hfill
    \begin{subfigure}[b]{0.24\textwidth}
        \centering
        \includegraphics[width=\textwidth]{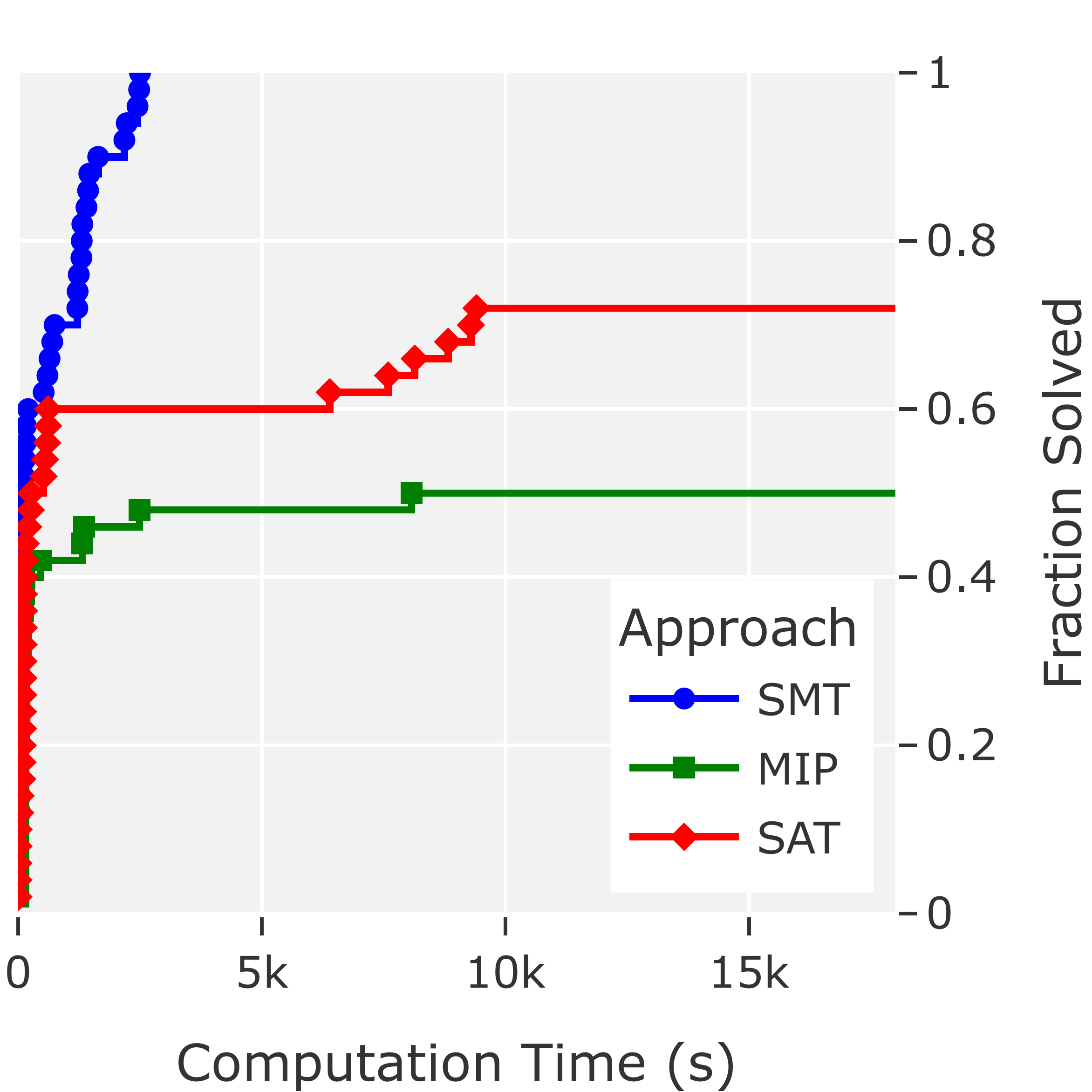}
        \caption{3-second granularity}
        \label{fig:cumulative_success_dt3_fs}
    \end{subfigure}
    \hfill
    \begin{subfigure}[b]{0.24\textwidth}
        \centering
        \includegraphics[width=\textwidth]{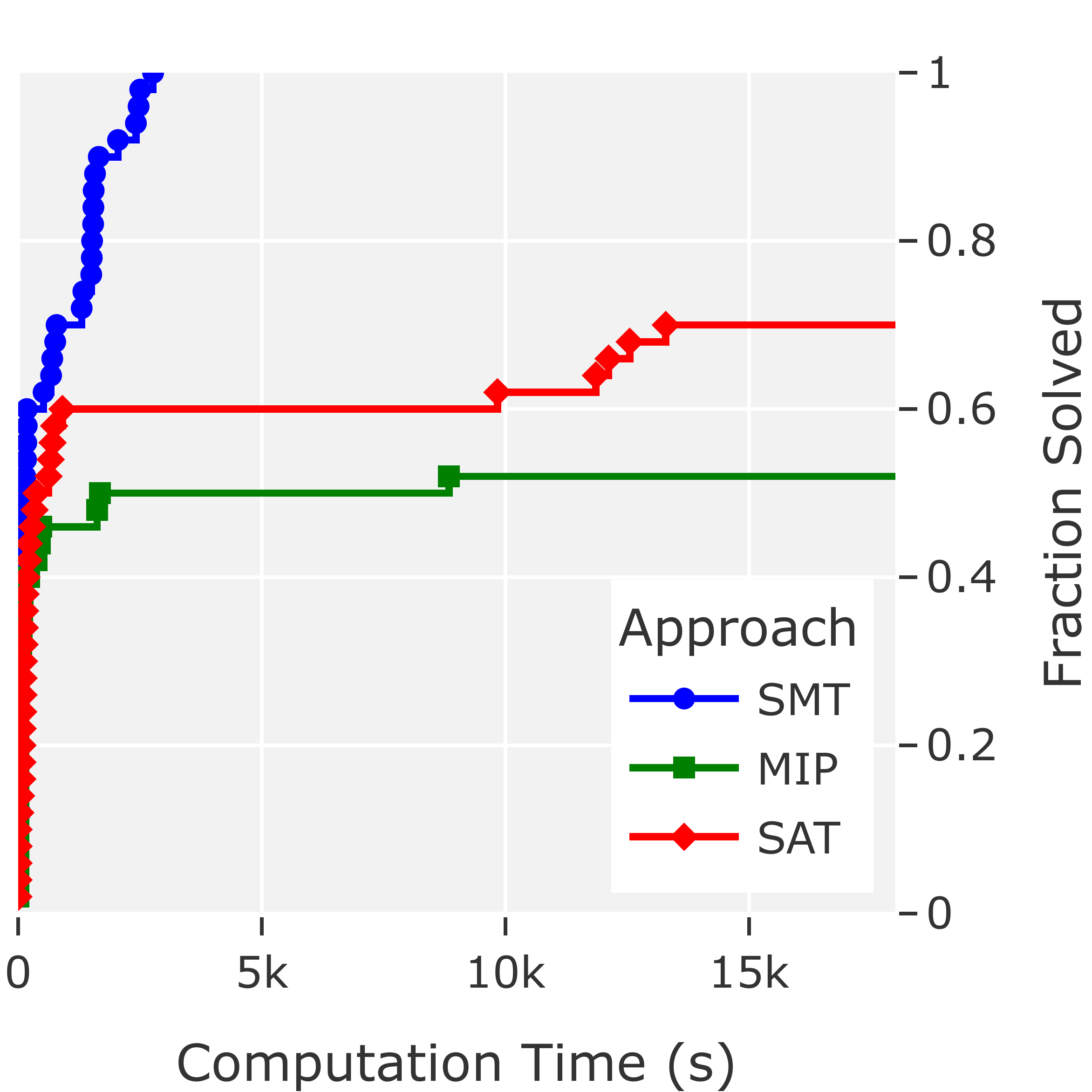}
        \caption{2-second granularity}
        \label{fig:cumulative_success_dt2_fs}
    \end{subfigure}
    \hfill
    \begin{subfigure}[b]{0.24\textwidth}
        \centering
        \includegraphics[width=\textwidth]{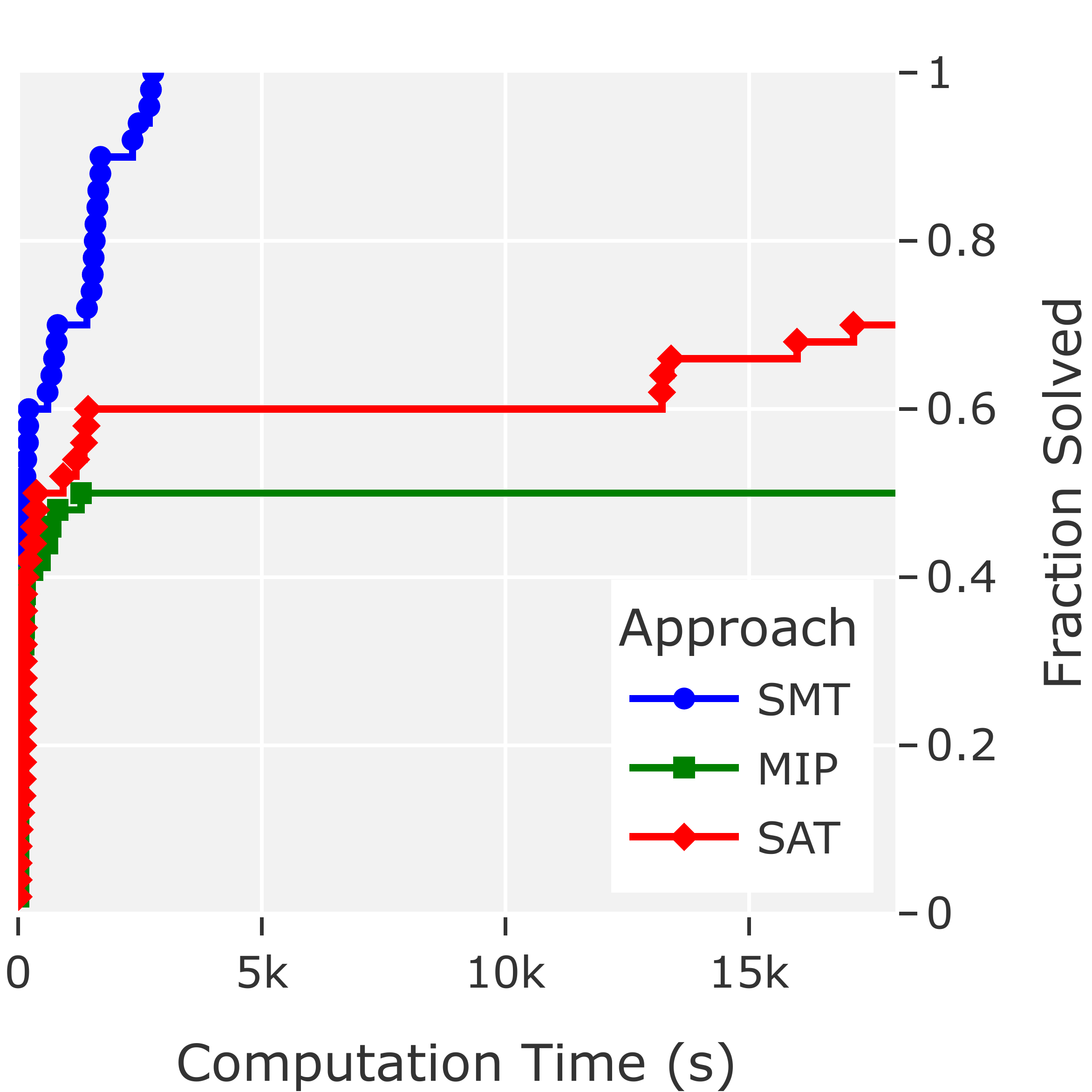}
        \caption{1-second granularity}
        \label{fig:cumulative_success_dt1_fs}
    \end{subfigure}
    \caption{Cumulative success rate for \acrshort{mip}, \acrshort{sat}, and \acrshort{smt} under \texttt{Full-Sharing} strategy.}
    \label{fig:cumulative_success_fs}
\end{figure}

The cumulative success plots in \cref{fig:cumulative_success_ns} (\texttt{No-Sharing}) and \cref{fig:cumulative_success_fs} (\texttt{Full-Sharing}) provide a more detailed view of solver performance across instances. \acrshort{smt} consistently solves all configurations within two hours, often well before the limit. \acrshort{sat} performs well at coarse granularities but exhibits a clear tipping point, beyond which runtime increases steeply and success rate drops. \acrshort{mip} solves the fewest instances overall, typically succeeding quickly or not at all.

The relative impact of sharing strategies is also visible. Under \texttt{Full-Sharing}, solution space complexity increases, and solvers generally take longer, particularly for \acrshort{sat}. Conversely, \acrshort{smt} 's completion profile remains identical, mainly, confirming its suitability for large and complex configurations.

These observations are consistent with the trends seen in computation times and are further reinforced by the success rate analysis. Across both circulation strategies, \acrshort{smt} emerges as the most reliable and scalable solver, maintaining low runtime and solving all instances. \acrshort{sat} performs well in simpler settings but scales poorly. \acrshort{mip} is only effective on small networks.

Time discretisation and vehicle sharing influence solver performance, with finer time steps and full-sharing policies introducing additional computational load.

\subsection{Effect of Time Discretisation and Vehicle Sharing Strategies on Vehicle Count}\label{ss:effect_discretisation}

This section analyses the impact of time discretisation on vehicle requirements, comparing the five different granularities used in Section~\ref{ss:performance_comparison}. While the baseline resolution is without discretisation, we apply four different levels of discretisation (1, 2, 3, and 6~seconds) to assess the impact of discretisation on the number of vehicles needed. Furthermore, we evaluate the potential benefits of increased flexibility through vehicle sharing by comparing both vehicle circulation strategies (\texttt{No-Sharing} and \texttt{Full-Sharing}).

When solving each scenario, we distinguish between two conflict configurations to isolate the effects of time discretization. This distinction helps separate the impact of rounding durations from the additional restrictions imposed by infrastructure conflicts. In real-world applications, rounding up activity durations is required for operational safety, but it can reduce schedule flexibility and increase travel times, potentially requiring more vehicles.

We therefore analyse:
\begin{enumerate}
    \item \texttt{No Conflicts}: Only activity durations (e.g., running and dwell times) are rounded up to match the discretisation step, while headway and other conflict-related constraints are omitted. This setup quantifies the isolated impact of rounding.
    \item \texttt{Conflicts}: All activities, including those modeling headways and resource usage, are adjusted for discretisation and fully included in the model. This setup reflects the full impact of time discretisation under realistic operational constraints.
\end{enumerate}

\Cref{tab:vehicle_counts} reports the required vehicles for each configuration. Parentheses indicate \texttt{No Conflicts} cases. We include the relative percentage increase in square brackets compared to the corresponding non-discretized baseline.

\begin{table}[H]
\centering
\caption{Vehicle requirements under varying time discretisation steps for increasing network sizes (given by the number of lines from 1 to 10). Values show absolute fleet size and the relative increase over the non-discretised baseline.}
\label{tab:vehicle_counts}
\begin{footnotesize}
\resizebox{\textwidth}{!}{%
\begin{tabular}{l|c|cccccccccc}
\toprule
\textbf{Strategy} & \textbf{Step} & \textbf{1} & \textbf{2} & \textbf{3} & \textbf{4} & \textbf{5} & \textbf{6} & \textbf{7} & \textbf{8} & \textbf{9} & \textbf{10} \\
\midrule
\multirow{5}{*}{\shortstack[l]{ 
\texttt{No-Sharing} \\ No Conflicts}}
& -- & 2 (0\%) & 4 (0\%) & 5 (0\%) & 7 (0\%) & 10 (0\%) & 14 (0\%) & 17 (0\%) & 20 (0\%) & 21 (0\%) & 24 (0\%) \\
& 1 sec  & 2 (0\%) & 4 (0\%) & 5 (0\%) & 7 (0\%) & 10 (0\%) & 14 (0\%) & 17 (0\%) & 20 (0\%) & 21 (0\%) & 24 (0\%) \\
& 2 sec  & 2 (0\%) & 4 (0\%) & 5 (0\%) & 7 (0\%) & 10 (0\%) & 14 (0\%) & 17 (0\%) & 20 (0\%) & 21 (0\%) & 24 (0\%) \\
& 3 sec  & 3 (+50.0\%) & 5 (+25.0\%) & 6 (+20.0\%) & 9 (+28.6\%) & 12 (+20.0\%) & 16 (+14.3\%) & 19 (+11.8\%) & 22 (+10.0\%) & 24 (+14.3\%) & 27 (+12.5\%) \\
& 6 sec  & 3 (+50.0\%) & 5 (+25.0\%) & 6 (+20.0\%) & 9 (+28.6\%) & 12 (+20.0\%) & 16 (+14.3\%) & 19 (+11.8\%) & 22 (+10.0\%) & 24 (+14.3\%) & 27 (+12.5\%) \\
\midrule
\multirow{5}{*}{
\shortstack[l]{ 
\texttt{No-Sharing} \\ \texttt{Conflicts}}}
& -- & 3 (0\%) & 5 (0\%) & 6 (0\%) & 9 (0\%) & 12 (0\%) & 16 (0\%) & 19 (0\%) & 22 (0\%) & 23 (0\%) & 26 (0\%) \\
& 1 sec  & 3 (0\%) & 5 (0\%) & 6 (0\%) & 9 (0\%) & 12 (0\%) & 16 (0\%) & 19 (0\%) & 22 (0\%) & 23 (0\%) & 26 (0\%) \\
& 2 sec  & 3 (0\%) & 5 (0\%) & 6 (0\%) & 9 (0\%) & 12 (0\%) & 16 (0\%) & 19 (0\%) & 22 (0\%) & 23 (0\%) & 26 (0\%) \\
& 3 sec  & 3 (0\%) & 5 (0\%) & 6 (0\%) & 9 (0\%) & 12 (0\%) & 16 (0\%) & 19 (0\%) & 22 (0\%) & 24 (+4.3\%) & 27 (+3.8\%) \\
& 6 sec  & 3 (0\%) & 5 (0\%) & 6 (0\%) & 9 (0\%) & 12 (0\%) & 17 (+6.3\%) & 20 (+5.3\%) & 23 (+4.5\%) & 25 (+8.7\%) & 28 (+7.7\%) \\
\midrule
\multirow{5}{*}{\shortstack[l]{ 
\texttt{Full-Sharing} \\ \texttt{No-Conflicts}}}
& -- & 2 (0\%) & 4 (0\%) & 4 (0\%) & 6 (0\%) & 9 (0\%) & 12 (0\%) & 15 (0\%) & 17 (0\%) & 18 (0\%) & 21 (0\%) \\
& 1 sec  & 2 (0\%) & 4 (0\%) & 4 (0\%) & 6 (0\%) & 9 (0\%) & 13 (+8.3\%) & 15 (0\%) & 17 (0\%) & 18 (0\%) & 21 (0\%) \\
& 2 sec  & 2 (0\%) & 4 (0\%) & 4 (0\%) & 6 (0\%) & 9 (0\%) & 13 (+8.3\%) & 15 (0\%) & 18 (+5.9\%) & 19 (+5.6\%) & 21 (0\%) \\
& 3 sec  & 3 (+50.0\%) & 4 (0\%) & 4 (0\%) & 7 (+16.7\%) & 10 (+11.1\%) & 14 (+16.7\%) & 16 (+6.7\%) & 19 (+11.8\%) & 20 (+11.1\%) & 22 (+4.8\%) \\
& 6 sec  & 3 (+50.0\%) & 4 (0\%) & 5 (+25.0\%) & 8 (+33.3\%) & 10 (+11.1\%) & 14 (+16.7\%) & 16 (+6.7\%) & 19 (+11.8\%) & 20 (+11.1\%) & 23 (+9.5\%) \\
\midrule
\multirow{5}{*}{\shortstack[l]{ 
\texttt{Full-Sharing} \\ \texttt{Conflicts}}}
& -- & 3 (0\%) & 4 (0\%) & 4 (0\%) & 7 (0\%) & 10 (0\%) & 14 (0\%) & 17 (0\%) & 20 (0\%) & 21 (0\%) & 24 (0\%) \\
& 1 sec  & 3 (0\%) & 4 (0\%) & 4 (0\%) & 7 (0\%) & 10 (0\%) & 14 (0\%) & 17 (0\%) & 20 (0\%) & 21 (0\%) & 24 (0\%) \\
& 2 sec  & 3 (0\%) & 4 (0\%) & 4 (0\%) & 7 (0\%) & 10 (0\%) & 14 (0\%) & 17 (0\%) & 20 (0\%) & 21 (0\%) & 24 (0\%) \\
& 3 sec  & 3 (0\%) & 4 (0\%) & 4 (0\%) & 7 (0\%) & 10 (0\%) & 14 (0\%) & 17 (0\%) & 20 (0\%) & 21 (0\%) & 24 (0\%) \\
& 6 sec  & 3 (0\%) & 4 (0\%) & 5 (+25.0\%) & 8 (+14.3\%) & 11 (+10.0\%) & 15 (+7.1\%) & 18 (+5.9\%) & 21 (+5.0\%) & 22 (+4.8\%) & 25 (+4.2\%) \\
\bottomrule
\end{tabular}
}
\end{footnotesize}
\end{table}

The results in \cref{tab:vehicle_counts} confirm several previous trends. The \texttt{Full-Sharing} strategy consistently results in lower fleet requirements than \texttt{No-Sharing}, as the additional flexibility enables more efficient vehicle reuse. This effect is evident across all line counts and discretisation levels. As expected, vehicle requirements increase with network size due to higher service numbers and more interactions due to shared infrastructure. On average, across all configurations, adopting \texttt{Full-Sharing} instead of \texttt{No-Sharing} reduces the vehicle count by 12.9\% in the \texttt{No Conflicts} setting and by 12.1\% when headway constraints are included.

Discretisation has a visible and cumulative impact. Under the \texttt{No Conflicts} setting, coarser time steps lead to systematic increases in fleet size. For example, under \texttt{No-Sharing}, shifting from the continuous-time baseline to 6-second discretisation increases total vehicle needs by 19 vehicles, corresponding to an average increase of 15.3\%. Under \texttt{Full-Sharing}, the increase is 14 vehicles or 13.0\%. These differences arise solely from reduced scheduling precision, without any infrastructure constraints.

When headway constraints are introduced, the relative increases are smaller: seven vehicles (+5.0\%) for \texttt{No-Sharing} and eight vehicles (+6.5\%) for \texttt{Full-Sharing}. In this case, constraints already restrict the solution space, making it less sensitive to discretisation. Nevertheless, finer time steps remain beneficial in minimizing the required fleet.

Notably, switching the vehicle circulation strategy cannot compensate for the increase due to coarser steps. Although \texttt{Full-Sharing} reduces vehicle needs overall, it does not eliminate the structural overhead introduced by time rounding. The cost of discretisation is additive to the cost of constraint-induced rigidity.

These findings underscore the importance of modeling both precise timing and flexible circulation. Coarse discretisation inflates vehicle requirements significantly, even without conflicts. The impact is somewhat dampened but still measurable when infrastructure constraints are present. Overall, modeling accuracy and resource flexibility are essential for efficient railway operations.

\subsection{Impact of Fixed Routes on Vehicle Requirements}\label{ss:impact_fixed_routes}

This section evaluates the effect of fixing train routes on vehicle requirements. While previous sections assume routing flexibility, we explore a constrained variant where routes are pre-assigned. To ensure feasibility, the fixed routes are extracted from solutions previously computed with full routing flexibility. Specifically, we solve the original problem under flexible routing and enforce the obtained paths in a second run. This enables us to assess the cost of eliminating routing freedom while avoiding infeasibility.

The experiments are run using the \acrshort{smt}-based solver with a one-second discretisation step under \texttt{No-Sharing} and \texttt{Full-Sharing} strategies. \cref{tab:impact_fixed_routing} reports the number of vehicles required in the fixed-routing case and the percentage increase compared to the flexible-routing baseline.

\begin{table}[H]
\centering
\caption{Vehicle requirements with fixed train routing for increasing network sizes (1–10 lines). Values show absolute fleet size, with relative increase over the flexible-routing baseline.}
\label{tab:impact_fixed_routing}
\resizebox{0.98\textwidth}{!}{%
\begin{tabular}{l|cccccccccc}
\toprule
\textbf{Strategy} & \textbf{1} & \textbf{2} & \textbf{3} & \textbf{4} & \textbf{5} & \textbf{6} & \textbf{7} & \textbf{8} & \textbf{9} & \textbf{10} \\
\midrule
\texttt{No-Sharing}     & 3 (+0.0\%) & 5 (+0.0\%) & 6 (+0.0\%) & 9\,~ (+0.0\%) & 13\,~ (+8.3\%) & 17 (+6.3\%) & 19 (+0.0\%) & 22 (+0.0\%) & 24 (+4.3\%) & 27 (+3.8\%) \\
\texttt{Full-Sharing}   & 3 (+0.0\%) & 4 (+0.0\%) & 5 (+25.0\%) & 8 (+14.3\%) & 11 (+10.0\%) & 14 (+0.0\%) & 17 (+0.0\%) & 21 (+5.0\%) & 22 (+4.8\%) & 24 (+0.0\%) \\
\bottomrule
\end{tabular}
}
\end{table}

The results in \cref{tab:impact_fixed_routing} confirm that fixing routes generally increases the required fleet size. For the \texttt{No-Sharing} strategy, the average increase is 2.3\%, while for \texttt{Full-Sharing} it is slightly higher at 5.9\%. This difference highlights that, while full-sharing provides more efficient vehicle utilisation under flexible routing, it is also more sensitive to restrictions once flexibility is removed. Nevertheless, even under fixed routing, in most cases, the configuration of the \texttt{Full-Sharing} configuration still requires fewer vehicles overall than the \texttt{No-Sharing} baseline, confirming that the inherent benefits of vehicle sharing persist despite routing constraints.

The impact of fixed routing remains negligible for small networks, where route alternatives are limited and the feasible space is less constrained. However, from medium-sized networks onward, the effect becomes visible. In the \texttt{Full-Sharing} case, pre-assigned paths require up to 25\% more vehicles than their flexible counterparts (e.g., instances with three lines). Interestingly, increases are not uniformly distributed. While some configurations (e.g., instances with 6 or 10 lines) show no penalty, others, such as with 4 or 5 lines, exhibit sharp increases. This suggests that the availability of alternative paths plays a key role: where the flexible model can exploit routing options to compress turnaround times or vehicle transitions, the fixed model is more likely to induce schedule fragmentation and idle time.

Overall, the findings reinforce the role of routing flexibility in reducing rolling stock requirements. While fixing paths may be operationally convenient, it should be applied with care in larger systems, as it can eliminate optimization potential and increase operational cost through larger fleet needs. The variability in impact across network sizes suggests that routing decisions should be informed by network structure and available topological alternatives, rather than treated as static inputs.

\subsection{Benefit of Integration}\label{ss:benefit_of_integration}

This section investigates the benefit of jointly optimizing train routing, periodic timetabling, and vehicle circulation, rather than solving these components sequentially. Traditional planning pipelines typically approach these problems in isolation: routes are fixed first, then a timetable is computed—often optimizing travel times—and finally a feasible vehicle circulation is derived. While intuitive, this neglects interdependencies that can impact resource efficiency.

An assumption sometimes made is that minimizing travel times indirectly reduces fleet size. For example, \citet{Liebchen2004TheBeyond} report that minimizing passenger travel time on the Berlin Underground led to fewer required vehicles, despite vehicle count not being part of the objective. Our study offers another perspective by explicitly targeting vehicle minimization.

\citet{lieshout2021VehiclePesp} previously demonstrated the benefit of integrating periodic timetabling and vehicle circulation regarding travel time and operational costs. In contrast, our work focuses exclusively on the required fleet size and extends the integrated formulation to include routing decisions. This allows us to isolate the impact of full integration while avoiding the need to model passenger travel times directly.

To quantify the effect, we compare two approaches:
\begin{itemize}
    \item \textbf{Integrated Solve:} Uses our model, where routing, timetabling, and circulation are jointly optimized to minimize the number of vehicles.
    \item \textbf{Sequential Solve:} Routing is fixed first, then timetabling is solved with an objective of minimal travel time, followed by vehicle circulation. Multiple runs with different initial seeds ensure robustness of the results.
\end{itemize}

The results are shown in \cref{tab:benefit_of_integration_table}, which reports the required fleet size for the sequential solve and its relative increase over the integrated baseline.

\begin{table}[H]
\centering
\caption{Vehicle requirements under sequential optimization for increasing network sizes (1–10 lines). Values show absolute fleet size and relative increase over the integrated baseline.}
\label{tab:benefit_of_integration_table}
\resizebox{0.98\textwidth}{!}{%
\begin{tabular}{l|cccccccccc}
\toprule
\textbf{Strategy} & \textbf{1} & \textbf{2} & \textbf{3} & \textbf{4} & \textbf{5} & \textbf{6} & \textbf{7} & \textbf{8} & \textbf{9} & \textbf{10} \\
\midrule
\texttt{No-Sharing}     & 3 (+0.0\%) & 5 (+0.0\%) & 7 (+16.7\%) & 9 (+0.0\%) & 13 (+8.3\%) & 19 (+18.8\%) & 21 (+10.5\%) & 25 (+13.6\%) & 29 (+26.1\%) & 32 (+23.1\%) \\
\texttt{Full-Sharing}   & 3 (+0.0\%) & 4 (+0.0\%) & 6 (+50.0\%) & 8 (+14.3\%) & 11 (+10.0\%) & 19 (+35.7\%) & 21 (+17.6\%) & 26 (+30.0\%) & 25 (+19.0\%) & 30 (+25.0\%) \\
\bottomrule
\end{tabular}
}
\end{table}

The results show that the integrated approach consistently reduces fleet size compared to the sequential strategy. Under the \texttt{No-Sharing} policy, the average reduction is 11.7\%, increasing with the number of lines. For \texttt{Full-Sharing}, the improvement reaches up to 35.7\%, highlighting how routing flexibility and vehicle sharing jointly benefit from integration. Notably, savings persist even in smaller configurations, with a marked increase in larger networks, where decisions across planning stages interact more tightly.

The increased vehicle count in the sequential solve stems from misalignments between routing, scheduling, and circulation. Minimizing travel time can yield structurally inefficient timetables, such as tight or incompatible turns, requiring additional vehicles to preserve feasibility. In contrast, the integrated model internalizes such constraints early and identifies globally efficient schedules.

While sequential approaches are easier to implement and align with traditional planning practice, they leave significant optimization potential untapped. Our results reinforce that jointly solving routing, scheduling, and circulation is essential to minimize operational cost and fleet size, especially as system complexity grows.

\section{Conclusion}\label{sec:conclusion}

This paper presents a novel, integrated optimization model for periodic timetabling, train routing, and vehicle circulation in railway systems. We propose the VCR-PESP formulation and introduce the first Satisfiability Modulo Theories (\acrshort{smt})-based solution method tailored to this problem. In contrast to traditional approaches using Mixed Integer Programming (\acrshort{mip}) or Boolean Satisfiability (\acrshort{sat})~\citep{Gromann2016SatisfiabilityProblems, Kummling2015AComputation}, our \acrshort{smt} formulation enables continuous-time modeling and supports scalable solving without time discretisation~\citep{armando2004sat, LEUTWILER2022525}.

Using real-world data from the Swiss narrow-gauge network of \acrshort{rhb}, we conduct an extensive computational study to assess the impact of time discretisation, routing flexibility, and vehicle sharing strategies. Our results show that coarser discretisations significantly increase fleet requirements due to rounding, and this effect persists even in conflict-free configurations. While vehicle sharing across lines can mitigate some of this overhead, only continuous-time models—as enabled by \acrshort{smt}—consistently avoid these losses without compromising resolution.

In all tested scenarios, our \acrshort{smt} solver outperforms \acrshort{sat} and \acrshort{mip}, solving large-scale instances with higher precision and shorter runtimes. We also show that restricting routing or solving planning stages sequentially leads to unnecessary vehicle use. In contrast, jointly optimizing routing, timetabling, and vehicle circulation within a single model yields more efficient and implementable solutions—confirming insights from earlier integration studies~\citep{lieshout2021VehiclePesp, Goossens2006OnProblems} and extending them to include routing.

Future work could extend our model to incorporate passenger-centric objectives such as travel time and transfer quality~\citep{Polinder2021,Gattermann2016IntegratingApproach, Kroon2014FlexibleTimetabling}, or to address robustness and real-time re-optimization. Our approach lays a strong foundation for high-resolution, scalable, and fully integrated railway planning using modern constraint solving techniques.

\section*{Author Contributions}
\textbf{Florian Fuchs}: Conceptualization, Methodology, Software, Investigation, Formal analysis, Visualization, Writing – original draft. \textbf{Bernardo Martin-Iradi}: Conceptualization, Methodology, Writing – review \& editing. \textbf{Francesco Corman}: Supervision, Project administration, Writing – review \& editing.

\textbf{Declaration of Generative AI and AI-assisted technologies in the writing process}\\
Statement: During the preparation of this work the author(s) used ChatGPT in order to occasionally help with readability. After using this tool/service, the author(s) reviewed and edited the content as needed and take(s) full responsibility for the publication's content.

\end{spacing}
\vspace{3mm}
\bibliographystyle{abbrvnat}
\bibliography{main.bbl}

\newacronym{los}{\texttt{LOS}}{Level of Service}

\newacronym{sq}{\texttt{SQ}}{Sub-Question (of the Research Question)}
\newacronym{rq}{\texttt{RQ}}{(Main) Research Question}
\newacronym{wp}{\texttt{WP}}{Work Package}

\newacronym{lbbd}{\texttt{LBBD}}{Logic-Based Benders Decomposition}
\newacronym{vcp}{\texttt{VCP}}{Vehicle Circulation Problem}

\newacronym{lpp}{\texttt{LPP}}{Line Planning Problem}
\newacronym{ttp}{\texttt{TTP}}{Train Timetabling Problem}
\newacronym{vsp}{\texttt{VSP}}{Vehicle Scheduling Problem}
\newacronym{sssp}{\texttt{SSSP}}{Slot Sequence Selection Problem}
\newacronym{m-sssp}{\texttt{M-SSSP}}{Master Slot Sequence Selection Problem}
\newacronym{s-sssp}{\texttt{S-SSSP}}{Sub Slot Sequence Selection Problem}
\newacronym{primal-sssp}{\texttt{Primal-SSSP}}{Primal Slot Sequence Selection Problem}
\newacronym{dual-sssp}{\texttt{Dual-SSSP}}{Dual Slot Sequence Selection Problem}

\newacronym{pesp}{\texttt{PESP}}{Periodic Event Scheduling Problem}
\newacronym{mip}{\texttt{MIP}}{Mixed Integer Program}
\newacronym{sat}{\texttt{SAT}}{Boolean Satisfiability Problem}
\newacronym{smt}{\texttt{SMT}}{Satisfiability Modulo Theories}

\newacronym{cgn}{\texttt{CGN}}{Change \& Go Network}
\newacronym{ptn}{\texttt{PTN}}{Public Transport Network}
\newacronym{cnf}{\texttt{CNF}}{Conjunctive Normal Form}
\newacronym{dag}{\texttt{DAG}}{Directed Acyclic Graph}

\newacronym{tpe}{\texttt{TPE}}{Train Path Envelope}
\newacronym{tss}{\texttt{TSS}}{Train Slot Sequence}
\newacronym{am}{\texttt{AM}}{Assessment Module}
\newacronym{om}{\texttt{OM}}{Optimisation Module}

\newacronym{toc}{\texttt{TOC}}{Train Operating Company}
\newacronym{im}{\texttt{IM}}{Infrastructure Manager}

\newacronym{rhb}{\texttt{RhB}}{Rhaetian Railway}

\newacronym{ean}{EAN}{Event Activity Network}
\newacronym{tfn}{TFN}{Train Flow Network}

\printglossary[type=\acronymtype]  

\newpage
\appendix

\end{document}